\documentclass[11pt]{amsart}
\usepackage{amsmath}
\usepackage{amssymb}

\usepackage{graphicx}
\usepackage[hidelinks]{hyperref}
\newcommand{\alg}[1]{{\boldsymbol{#1}}} % algebras
\newcommand{\cl}[1]{{\mathsf{#1}}} % classes

\newcommand{\leibniz}{\boldsymbol{\varOmega}}

\newcommand{\Con}{\operatorname{Con}}
\newcommand{\ConL}{\boldsymbol{{\operatorname{Con}}}}
\newcommand{\DFil}{\operatorname{Fil}}
\newcommand{\DFilL}{\boldsymbol{{\operatorname{Fil}}}}

\newcommand{\Sg}{\operatorname{Sg}}
\newcommand{\SgA}{\boldsymbol{{\operatorname{Sg}}}}

\newcommand{\Vop}{\mathbb{V}}

\newcommand{\Hop}{\mathbb{H}}
\newcommand{\Sop}{\mathbb{S}}

\newcommand{\Puop}{\mathbb{P}_\mathbb{U}}
\newcommand{\Psop}{\mathbb{P}_\mathbb{S}}
\newcommand{\Iop}{\mathbb{I}}
\newcommand{\Rop}{\mathbb{R}}

\newcommand{\refl}{\operatorname{R}}
\newcommand{\bdot}{\mathbin{\boldsymbol{\cdot}}}
\newcommand{\seteq}{\mathrel{\mbox{\,\textup{:}\!}=\nolinebreak }\,}

\newcommand{\0}{\mathsf{\mathit{0}}}
\newcommand{\blank}{\,\textup{-}\,}

\newcommand{\id}{{\textup{id}}}
\newcommand{\aand}{\mathbin{\;\&\;}}
\newcommand{\Z}{\mathbb{Z}}

\newtheorem{theorem}{Theorem}[section] %TODO change to [chapter] in thesis
\newtheorem{lemma}[theorem]{Lemma}
\newtheorem{corollary}[theorem]{Corollary}

\theoremstyle{definition}
\newtheorem{definition}[theorem]{Definition}

\newtheorem{example}[theorem]{Example}

\makeatletter

\renewcommand{\theenumi}{\roman{enumi}}
\renewcommand{\theenumii}{\roman{enumii}}

\renewcommand{\p@enumii}{\theenumi(\theenumii)}
\makeatother

\title[Semilinear De Morgan monoids and epimorphisms]{Semilinear De Morgan monoids and epimorphisms}
\author{J.J.\ Wannenburg}
\address{%
Institute of Computer Science, Academy of Sciences of the Czech Republic, Pod Vod\'{a}renskou v\v{e}\v{z}\'{i} 2, 182 07 Prague 8, Czech Republic
}
\email{wannenburg@cs.cas.cz}
\author{J.G.\ Raftery}
\address{Department of Mathematics and Applied Mathematics,
University of Pretoria,
Private Bag X20, Hatfield,
Pretoria 0028, South Africa}
\email{{james.raftery@up.ac.za}} % \quad {jamie.wannenburg@up.ac.za}
\thanks{%
The first author's work was carried out within the project \emph{Supporting the
internationalization of the Institute of Computer Science of the Czech
Academy of Sciences} (no.~CZ.02.2.69/0.0/0.0/18\_053/0017594), funded by
the Operational Programme Research, Development and Education of the
Ministry of Education, Youth and Sports of the Czech Republic. The
project is co-funded by the EU.  The second author was supported in part
by the National Research Foundation of South Africa (UID 85407).
}

\begin{document}

\begin{abstract}
A representation theorem is
proved for 
De Morgan monoids 
that are
(i) \emph{semilinear}, i.e., subdirect products of totally ordered algebras, and (ii) \emph{negatively generated}, i.e., generated by lower bounds of the neutral element.  Using this theorem, we prove that the De Morgan monoids satisfying (i) and (ii) form a variety---in fact, a locally finite variety.  We then prove that epimorphisms are surjective in every variety of negatively generated semilinear De Morgan monoids.  In the process, epimorphism-surjectivity is established for several other classes as well, including
the variety of all semilinear idempotent commutative residuated
lattices and all varieties of negatively generated semilinear Dunn monoids.
The results settle natural questions about Beth-style definability for a range of substructural logics. 
\end{abstract}

\maketitle

%\makeatletter
%\renewcommand{\labelenumi}{\text{(\theenumi)}}
%\renewcommand{\theenumi}{\roman{enumi}}
%\renewcommand{\theenumii}{\roman{enumii}}
%\renewcommand{\labelenumii}{\text{(\theenumii)}}
%\renewcommand{\p@enumii}{\theenumi(\theenumii)}
%\makeatother

%\allowdisplaybreaks

\section{Introduction}
\label{sec:introduction}
The aims of this paper are two-fold.  On one hand, we continue a line of investigation in \cite{BMR17,MRW4,MW}, which seeks to identify varieties of residuated structures in which epimorphisms are surjective (a property that need not persist in subvarieties).  To do this, we need to prove some structural representation theorems for the algebras concerned.
These representation theorems are of independent interest, and their exposure is our second aim.

The first aim is motivated by the connection between Beth-style definability properties in substructural logics and the behaviour of epimorphisms in the varieties of residuated structures that model them.  We are concerned here with the surjectivity of \emph{all} epimorphisms in these varieties.  This characterizes the so-called \emph{infinite Beth property} for the corresponding logics \cite{BH06}; see \cite[pp.\,186--7]{BMR17} for a concise account of the details, as well as references.

We are also concerned here with residuated structures that need not be \emph{integral}, i.e., the neutral element $e$ for fusion ($\bdot$) need not be the greatest element of the algebra.  Earlier investigations focussed mainly on strong or weak variants of epimorphism-surjectivity, and on integral structures, such as Heyting or Brouwerian algebras, which model intuitionistic propositional logic and its positive fragment \cite{EG81,Kre60,Mak03}.  The present study, like \cite{BMR17,MRW4}, accommodates various extensions of relevance logic as well, so the structures under consideration will be De Morgan and Dunn monoids.  (A De Morgan monoid is essentially a Dunn monoid equipped with an involution $\neg$ that simulates negation.)

It was proved in \cite[Thm.~8.1]{MRW4} that epimorphisms will be surjective in a variety of De Morgan or Dunn monoids, provided that the finitely subdirectly irreducible members of the variety are \emph{negatively generated} (i.e., generated by lower bounds of $e$) and that their posets of prime filters have finite depth.  We show here that the demand for finite depth can be dropped when the algebras are \emph{semilinear} (i.e., subdirect products of chains).

Whereas De Morgan and Dunn monoids satisfy the \emph{square-increasing law} $x\leqslant x^2\seteq x\bdot x$,
the negatively generated semilinear Dunn monoids turn out to be \emph{idempotent} (Theorem~\ref{thm:SLDunn-neg-cone}), i.e., they satisfy $x=x^2$.
They therefore coincide with the
\emph{generalized Sugihara monoids} of \cite{GR15,BMR17}, which form a locally finite variety $\cl{GSM}$.
We show that \emph{all} subvarieties of $\cl{GSM}$ have surjective epimorphisms (Theorem~\ref{thm:GSM-ES}).  We also show in Theorem~\ref{thm:SLid-ES}
that epimorphisms are
surjective in the variety of \emph{all} semilinear idempotent Dunn monoids (regardless of negative generation).  As that variety is known to have the amalgamation property \cite{GJM19}, it follows that it has the strong amalgamation property (Corollary~\ref{cor:indemSemiCRLhasStrongESandStrongAmalgamation}).

Using some of these results and a characterization of irreducible De Morgan monoids from \cite{MRW}, we show in Section~\ref{sec:DMM} that the negatively generated semilinear De Morgan monoids also form a variety---in fact a locally finite one (Corollary~\ref{cor:SL-DMM-loc-fin}).  We conclude with a proof that epimorphisms are surjective in each of its subvarieties (Theorem~\ref{thm:SL-DMM-ES}).

The key to this proof
is a structural result, Theorem~\ref{thm:SL-DMM-neg-cone}.  It says that every negatively generated and totally ordered De Morgan monoid $\alg{A}$ arises from a totally ordered generalized Sugihara monoid $\alg{G}$ by the use of two constructions.  First, we construct a `reflection' $\refl(\alg{G})$ of $\alg{G}$ which places an inverted copy of $\alg{G}$ above all the elements of $\alg{G}$, adding bounds and an involution, and extending the original operations systematically.  Secondly, we `substitute' $\refl(\alg{G})$ (in a suitable sense) for the neutral element of a totally ordered `odd Sugihara monoid' $\alg{S}$ (i.e., an idempotent De Morgan monoid in which $e=\neg e$).  We call this second construction a `rigorous extension'.  Because totally ordered odd Sugihara monoids are transparently structured, the resulting algebra $\alg{A}$ is easily analysed.  Relative to $\alg{A}$, the algebras $\refl(\alg{G})$ and $\alg{S}$ are the interval $[\neg((\neg e)^2),(\neg e)^2]$ and
the factor algebra
got by collapsing $\refl(\alg{G})$ to a point and isolating all other elements.

\section{Conventions}\label{sec:conventions}

As usual, $\omega$ denotes the set of non-negative integers.
The universe of an algebra $\alg{A}$ is denoted by $A$.  Thus, the congruence lattice $\ConL\,\alg{A}$
of $\alg{A}$ has universe $\Con\alg{A}$.
For $\emptyset\neq X\subseteq A$,
the subalgebra of $\alg{A}$ generated by $X$ is
denoted by $\SgA^\alg{A}X$
(and its universe by
$\Sg^\alg{A} X$).  An algebra $\alg{A}$ is said to be $n$-\emph{generated}, where $n\in\omega$, if it has the form $\SgA^\alg{A}X$ for some $X$ such that $|X|\leqslant n$.

The class operator symbols
$\Iop$, $\mathbb{H}$, $\mathbb{S}$, $\mathbb{P}$, $\Psop$
and $\mathbb{P}_\mathbb{U}$
stand, respectively, for closure under isomorphic and
homomorphic images, subalgebras,
direct and subdirect
products,
and ultra\-products,
while $\mathbb{V}$
denotes varietal
generation, i.e., $\mathbb{V}=\mathbb{HSP}$.
We abbreviate $\mathbb{V}(\{\alg{A}\})$ as $\mathbb{V}(\alg{A})$.

Recall that an algebra $\alg{A}$ is \emph{subdirectly irreducible} (SI)
iff its identity relation
$\id_A=\{\langle a,a\rangle : a\in A\}$
is completely
meet-irreducible in its congruence lattice.
Also, $\alg{A}$ is \emph{finitely subdirectly irreducible} (FSI)
iff\/ $\id_A$ is meet-irreducible in
$\ConL\,\alg{A}$,
whereas $\alg{A}$ is \emph{simple}
iff\/ $\left|\Con\alg{A}\right|=2$.  Consequently, trivial algebras are FSI, but are neither SI nor simple.

Let $\cl{K}$ be a variety.
We denote by
$\cl{K}_\textup{SI}$ [resp.\ $\cl{K}_\textup{FSI}$]
the class of subdirectly irreducible [resp.\ finitely subdirectly irreducible]
members of $\cl{K}$.  Thus, $\cl{K}=\mathbb{V}(\cl{K}_\textup{SI})$.
\emph{J\'{o}nsson's Theorem} \cite{Jon67,Jon95}
states that,
for any subclass
$\cl{L}$ of a
congruence distributive variety,
$\Vop(\cl{L})_\textup{FSI}
\subseteq
\Hop\Sop\Puop(\cl{L})$.
In this connection, recall that $\Puop(\cl{L})\subseteq\Iop(\cl{L})$ whenever $\cl{L}$ is a finite set of finite similar algebras.

\section{Epimorphisms}

Given a class $\cl{K}$ of similar algebras, a \emph{$\cl{K}$-morphism} is a homomorphism $\textup{$f \colon \alg{A} \to \alg{B}$}$,
where $\alg{A},\alg{B} \in \cl{K}$.  It is called a
\emph{$\cl{K}$-epimorphism} provided that, whenever
$g,h\colon\alg{B}\to\alg{C}$ are $\cl{K}$-morphisms
with $g \circ f = h \circ f$, then $g = h$.
Clearly, surjective $\cl{K}$-morphisms are $\cl{K}$-epimorphisms.  We say that
$\cl{K}$ has the \emph{epimorphism-surjectivity (ES) property} if all $\cl{K}$-epimorphisms are surjective.

A subalgebra $\alg{D}$ of an algebra $\alg{E}\in\cl{K}$ is said to be \emph{$\cl{K}$-epic} (in $\alg{E}$)
if every $\cl{K}$-morphism with
domain $\alg{E}$ is determined by its restriction to $\alg{D}$.  (This means that the inclusion map $\alg{D}\to\alg{E}$
is a $\cl{K}$-epimorphism, assuming that $\alg{D}\in\cl{K}$.)
Thus,
a $\cl{K}$-morphism is a $\cl{K}$-epimorphism iff its image is a $\cl{K}$-epic subalgebra of its co-domain.
And, when $\cl{K}$ is closed under subalgebras (in particular, when $\cl{K}$ is a variety),
then
\[
\emph{$\cl{K}$ has the ES property iff its members all lack $\cl{K}$-epic proper subalgebras.}
\]

A variety $\cl{K}$ is said to have \emph{EDPM} if it is congruence distributive and $\cl{K}_\textup{FSI}$ is a universal class (i.e.,
subalgebras and ultraproducts of FSI members of $\cl{K}$ are FSI).  The acronym stands for `equationally definable principal meets'
and is motivated by other characterizations of the notion in \cite{BP86,CD90}.

\begin{theorem}[{Campercholi \cite[Thm.~6.8]{Cam18}}]
\label{thm:Campercholi}
If a congruence permutable variety\/ $\cl{K}$ with EDPM lacks the ES property, then
some FSI member of
$\cl{K}$
has a\/ $\cl{K}$-epic proper subalgebra.
\end{theorem}

When testing whether a subalgebra is epic, we may also use the following consequence of the Subdirect Decomposition Theorem.

\begin{lemma}
\label{lem:FSI-image}
Let $\cl{K}$ be a variety of algebras and let $\alg{B}$ be a subalgebra of $\alg{A} \in \cl{K}$. Then $\alg{B}$ is\/ $\cl{K}$-epic in $\alg{A}$ iff, whenever $\alg{C} \in \cl{K}_{\textup{SI}}$ and $g,h \colon \alg{A} \rightarrow \alg{C}$ are homomorphisms that agree on $B$, then $g = h$.
\end{lemma}

\begin{definition}
\label{def:other-properties} Let $\cl{K}$ be a class of similar algebras.
\begin{enumerate}
\item \label{def:other-properties:weak-ES}
We say that $\cl{K}$ has the \emph{weak ES property} if
no finitely generated member of $\cl{K}$ has a $\cl{K}$-epic proper subalgebra.  An equivalent demand is that no $\alg{B}\in\cl{K}$ has a $\cl{K}$-epic proper subalgebra $\alg{A}$ such that
$B=\Sg^\alg{B}(A\cup C)$ for some \emph{finite} $C\subseteq B$ \cite[Thm.\,5.4]{MRW3}.

\smallskip

\item \label{def:other-properties:strong-ES} The \emph{strong ES property}
for $\cl{K}$ asks that, whenever $\alg{A}$ is a subalgebra of $\alg{B} \in \cl{K}$ and $b \in B \setminus A$, then there exist $\alg{C} \in \cl{K}$ and homomorphisms $g,h \colon \alg{B} \rightarrow \alg{C}$ such that $g|_{A} = h|_{A}$ and $g(b) \neq h(b)$.

\smallskip

\item \label{def:other-properties:amalgamation} The \emph{amalgamation property}
for a variety $\cl{K}$ is the demand that, for any two embeddings $g_B\colon \alg{A} \rightarrow \alg{B}$ and $g_C\colon \alg{A} \rightarrow \alg{C}$ between algebras in $\cl{K}$, there exist embeddings $h_B \colon \alg{B} \rightarrow \alg{D}$ and $h_C \colon \alg{C} \rightarrow \alg{D}$, with $\alg{D} \in \cl{K}$, such that $h_B \circ g_B = h_C \circ g_C$.

\smallskip

\item \label{def:other-properties:strong-amalgamation} The \emph{strong amalgamation property}
for $\cl{K}$ asks, in addition to the demands of (\ref{def:other-properties:amalgamation}), that $\alg{D}$, $h_B$ and $h_C$ can be chosen so that
\[
(h_B \circ g_B)[A] = h_B[B] \cap h_C[C].
\]
\end{enumerate}
\end{definition}
These conditions are linked as follows (see \cite{Isb66,Rin72,KMPT83} and \cite[Sec.~2.5.3]{Hoo01}).

\begin{theorem}
\label{thm:linkBetweenESandAmalgamation}
A variety has the strong amalgamation property iff it has the amalgamation property and the weak ES property. In that case, it also has the strong ES property (and therefore the ES property).
\end{theorem}

\section{Residuated Structures}
\label{sec:preliminaries}

\begin{definition}\label{def:IRL}
An \emph{involutive (commutative) residuated lattice}, or briefly, an \emph{IRL},
is an algebra $\alg{A}=\langle A;\bdot,\wedge,\vee,\neg,e\rangle$
comprising a commutative monoid $\langle A;\bdot,e\rangle$,
a lattice $\langle A;\wedge,\vee\rangle$
and a function $\neg\colon A\rightarrow A$,
called an \emph{involution}, 
such that $\alg{A}$ satisfies the (first order) formulas $\neg\neg x = x$ and
\begin{equation}\label{eq:neg-fusion-law}
x\bdot y\leqslant z\;\iff\;\neg z\bdot y\leqslant\neg x,
\end{equation}
cf.\ \cite{GJKO07}.\footnote{\,The signature in \cite{GJKO07} is slightly different, but the definable terms are not affected.}
Here, $\leqslant$ denotes the lattice order (i.e., $x\leqslant y$ abbreviates $x\wedge y= x$)
and $\neg$ binds more strongly than any other operation; we refer to $\bdot$
as \emph{fusion}.
\end{definition}
Setting $y=e$ in (\ref{eq:neg-fusion-law}), we see that $\neg$ is antitone.  In fact, De Morgan's laws for $\neg,\wedge,\vee$ hold, so
$\neg$ is an anti-automorphism of $\langle A;\wedge,\vee\rangle$.
If we define
\[
x\to y\seteq\neg(x\bdot\neg y)
\text{
 \ and \ } f\seteq\neg e,
\]
then, as is well known, every IRL satisfies
\begin{align}
& x\bdot y\leqslant z\;\Longleftrightarrow\;y\leqslant x\to z \quad\text{(the law of residuation),}\label{eq:residuation}
\\
& \neg x=  x\to f, \text{ \ hence \ } x\bdot\neg x\leqslant f,\label{eq:neg-via-f}\\
& x\to y= \neg y\to\neg x \text{ \ and \ } x\bdot y= \neg(x\to\neg y).\label{eq:contraposition-and-dot-via-arrow}
\end{align}

\begin{definition}\label{def:RL}
A \emph{(commutative) residuated lattice}---or
an \emph{RL}---is an algebra
$\alg{A}=\langle A;\bdot,\to,\wedge,\vee,e\rangle$
comprising a commutative monoid $\langle A;\bdot,e\rangle$, a lattice $\langle A;\wedge,\vee\rangle$ and
a binary operation $\to$, called the \emph{residual} of $\alg{A}$, 
where
$\alg{A}$ satisfies (\ref{eq:residuation}).
\end{definition}
Thus, up to term equivalence, every IRL has a reduct that is an RL. 
Conversely,
every RL can be embedded into (the RL-reduct of) an IRL; see \cite{GR04} and the
antecedents cited there.  Whereas $\bdot$ and $\to$ are inter-definable in IRLs, $\to$ is determined in RLs by $\bdot,\leqslant$,
because $x\to y$ coincides with ${\textup{max}\,\{z:x\bdot z\leqslant y\}}$.
Every RL satisfies the following well known formulas.
Here and subsequently, $x\leftrightarrow y$ abbreviates ${(x\to y)\wedge(y\to x)}$. \begin{align}
& x\bdot (x\to y)\leqslant y
\text{ \ and \ }  x\leqslant (x\to y)\to y
\label{eq:x-y-law}\\
& ((x \to y) \to y) \to y  =  x \to y \label{eq:triple-equality}\\
& (x\bdot y)\to z= y\to (x\to z)= x\to (y\to z) \label{eq:permutation}\\
& x\leqslant y\;\Longrightarrow\;
\begin{cases}
                           x\bdot z\leqslant y\bdot z
                           \;\text{ and }\;
                           \\
                           z\to x\leqslant z\to y
                           \;\text{ and }\;
                           y\to z\leqslant x\to z
\end{cases}
\label{eq:isotone}
\\
& x\leqslant y\;\Longleftrightarrow\;e\leqslant x\to y \label{eq:e-order}\\
& x= y\;\Longleftrightarrow\;e\leqslant x\leftrightarrow y \label{eq:e-equality}\\
& e\leqslant x\to x \text{ \ and \ } e\to x= x \label{eq:e-laws}.
\end{align}
The respective classes of all RLs and of all IRLs
are finitely axiomatizable varieties \cite[Thm.~2.7]{GJKO07}.

In an RL,
we define $x^0\seteq e$ and $x^{n+1}\seteq x^n\bdot x$ \,for $n\in\omega$.
\begin{definition}\label{def:S[I]RL}
An [I]RL is said to be \emph{square-increasing}
if it satisfies
\begin{equation}
\label{eq:square-increasing}
x\leqslant x^2 \quad \text{(the \emph{square-increasing law})}
\end{equation}
\end{definition}
Every square-increasing
RL can be embedded into a square-increasing IRL; see \cite{Mey73} and the `reflection' construction in Section~\ref{sec:DMM} below.
The following formulas are valid in all square-increasing RLs (and not in all RLs):
\begin{align}
& x\wedge y\leqslant x\bdot y
\label{eq:meet-below-fusion} \\
& (x\leqslant e \aand y\leqslant e)\;\Longrightarrow \; x\bdot y= x\wedge y.\label{eq:meet=fusion}
\end{align}

The next result is well known; see \cite[Cor.~14]{GOR08} and \cite[Thm.~2.4]{OR07}, for instance.

\begin{lemma}\label{lem:RL-properties}\
\begin{enumerate}
\item\label{lem:RL-properties:FSI}
An [I]RL
$\alg{A}$ is FSI iff\/ $e$ is join-irreducible in\/ $\langle A;\wedge,\vee\rangle$\textup{.}

\smallskip

\item\label{lem:RL-properties:SI}
A square-increasing [I]RL $\alg{A}$ is SI iff, in $\langle A;\wedge,\vee\rangle,$ there is a largest element strictly below $e$\textup.

\smallskip

\item\label{lem:RL-properties:simple}
A square-increasing [I]RL $\alg{A}$ is simple iff\/ $e$ has just one strict lower bound in $\langle A;\wedge,\vee\rangle$\textup.
\end{enumerate}
\end{lemma}

As RLs have lattice reducts, any
variety of
[I]RLs is congruence distributive.
It is also congruence permutable
and has the
congruence extension property (CEP);
see, for instance, \cite[Sec.~2.2 and 3.6]{GJKO07}.
Moreover, since the join-irreducibility of $e$ in condition~(\ref{lem:RL-properties:FSI}) is expressible as a universal first order sentence,
every variety
of
[I]RLs
has EDPM, so both J\'{o}nsson's Theorem and Theorem~\ref{thm:Campercholi} apply to such varieties.

An element $a$ of an [I]RL $\alg{A}$ is said to be \emph{idempotent}
if $a^2=a$.
We say that $\alg{A}$ is \emph{idempotent}
if all of its elements are.

Recall that for IRLs we defined the nullary term $f$ as $\neg e$. In \cite[Lem.~3.1]{MRW} it is shown that $f^3=f^2$ in any square-increasing IRL.  The following is stated in \cite[p.\,309]{Mey86}; for a proof see \cite[Lem.~3.3]{MRW}.

\begin{theorem}
\label{thm:idempotence}
In a square-increasing IRL $\alg{A}$\textup{,} the following are equivalent.
\begin{enumerate}
\item\label{thm:idempotence:f^2=f}
$f^2=f$\textup{.}
\item\label{thm:idempotence:f<=e}
$f\leqslant e$\textup{.}
\item\label{thm:idempotence:idempotent}
$\alg{A}$ is idempotent.
\end{enumerate}
Consequently, a square-increasing non-idempotent IRL has no idempotent subalgebra (and in particular, no trivial subalgebra).
\end{theorem}

An [I]RL $\alg{A}$ is said to be \emph{distributive} 
if its reduct $\langle A;\wedge,\vee\rangle$ is a
distributive 
lattice.
It is said to be
\emph{semilinear}
if it is
isomorphic to a subdirect product of totally ordered algebras (in which case it is obviously distributive).
Because the totality of a partial order is expressible by a universal positive sentence,
J\'{o}nsson's Theorem
has the following consequence:

\begin{lemma}
\label{lem:total-order}
A semilinear
[I]RL
$\alg{A}$ is FSI iff it is totally ordered.
\end{lemma}

It is shown in \cite{HRT02} that an [I]RL $\alg{A}$ is semilinear 
iff it is distributive and satisfies
$e\leqslant (x\to y)\vee(y\to x)$,
whence the semilinear [I]RLs form a variety.

Let $\alg{A}$ be an
[I]RL.
By a \emph{filter}
of
$\alg{A}$, we mean a filter of the lattice $\langle A;\wedge,\vee\rangle$, i.e.,
a non-empty subset $G$ of $A$ that is upward closed
and closed under the binary operation $\wedge$.
A \emph{deductive filter} 
of $\alg{A}$ is a filter $G$ of $\langle A;\wedge,\vee\rangle$ that is also a \emph{submonoid}
of $\langle A;\bdot,e\rangle$, i.e., $e \in G$ and $a \bdot b \in G$ whenever $a,b \in G$.
Thus, $[e) \seteq \{x \in A: x \geqslant e\}$ is the smallest deductive filter of $\alg{A}$,
and whenever $b\in A$ and $a,a\to b\in G$, then $b\in G$
(as $a\bdot(a\to b)\leqslant b$, by (\ref{eq:residuation})).
The lattice $\DFilL\,\alg{A}$
of deductive filters of $\alg{A}$ and the congruence lattice
$\ConL\,\alg{A}$
of $\alg{A}$ are isomorphic.  The isomorphism
and its inverse are given by
\begin{align*}
G\,& \mapsto\,\leibniz^{\alg{A}} G
\seteq
\{\langle a,b\rangle\in A^2 : a\to b,\,b\to a\in G\};\\
\theta\,& \mapsto\,\{a\in A : \langle a\wedge e, e\rangle\in\theta\}.
\end{align*}
For a deductive filter $G$ of $\alg{A}$ and $a,b\in A$, we often abbreviate $\alg{A}/\leibniz^{\alg{A}} G$ as $\alg{A}/G$,
and $a/\leibniz^{\alg{A}} G$ as $a/G$,
noting that
\begin{equation*}\label{eq:factor-order}
\textup{$a\to b\in G$ \,iff\,  $a/G\leqslant b/G$ in $\alg{A}/G$.}
\end{equation*}

When $\alg{A}$ is square-increasing,
the deductive filters of $\alg{A}$ are just the lattice filters of $\langle A;\wedge,\vee\rangle$ that contain $e$, by
(\ref{eq:meet-below-fusion}).

\section{De Morgan Monoids, Dunn Monoids and Sugihara Monoids}

\begin{definition}\label{def:DMM}
A \emph{De Morgan monoid}
is a distributive square-increasing IRL.
A \emph{Dunn monoid} is a distributive square-increasing RL.
\end{definition}

A \emph{Sugihara monoid}
is an idempotent De Morgan monoid, i.e.,
an idempotent distributive IRL.
The structure of such an algebra is better understood than that of an arbitrary De Morgan monoid,
largely
because of J.M.\ Dunn's
contributions
to \cite{AB75}; see \cite{Dun70} also.
The variety $\cl{SM}$ of all Sugihara monoids
is \emph{locally finite} (i.e., its finitely generated members are finite), but not \emph{finitely generated} (i.e., generated by a finite algebra).
In fact, $\cl{SM}$ is the smallest variety containing the Sugihara monoid
\[
\alg{Z}^*=\langle\{a:0\neq a\in\Z\};\bdot,\wedge,\vee,-,1\rangle
\]
on the set of all nonzero integers such that the lattice order is the usual total
order, the involution $-$ is the usual additive inversion, and the monoid operation is defined by
\[
a\bdot b \; = \; \begin{cases}
\text{the element of $\{a,b\}$ with the greater absolute value, if $\left|a\right| \neq \left|b\right|$;}\\
                           \text{$a\wedge b$  \,if $\left|a\right| = \left|b\right|$}
\end{cases}
\]
(where $\left|\blank\right|$ is the natural absolute value function). In this algebra, the residual operation $\to$ is given by
\begin{equation*}
a \to b \; = \; \begin{cases}
                           (-a) \vee b   & \text{if\, $a \leqslant b$;} \\
                           (-a) \wedge b & \text{if\, $a \not\leqslant b$.}
\end{cases}
\end{equation*}
Note that $e=1$ and $f=-1$ in $\alg{Z}^*$.

Because $\alg{Z}^{*}$ is totally ordered and generates $\cl{SM}$,
every FSI Sugihara monoid is totally ordered, i.e., Sugihara monoids are semilinear.

An IRL $\alg{A}$ is said to be \emph{odd}
if $f=e$ in $\alg{A}$.  Theorem~\ref{thm:idempotence} has the following consequence.

\begin{theorem}\label{thm:odd-DMMs-are-SM}
Every odd De Morgan monoid is a Sugihara monoid.
\end{theorem}

In the Sugihara monoid $\alg{Z}=\langle\Z;\bdot,\wedge,\vee,-,0\rangle$
on the set of \emph{all} integers, the operations are defined like those of $\alg{Z}^*$, except that
$0$ takes over from $1$ as the neutral element for $\bdot$.  Both $e$
and $f$ are $0$ in $\alg{Z}$, so $\alg{Z}$ is odd.
It follows from Theorem~\ref{thm:odd-DMMs-are-SM} and Dunn's results in \cite{AB75,Dun70} that the variety $\cl{OSM}$ of all odd Sugihara monoids
is the
smallest quasivariety containing $\alg{Z}$, and that $\cl{SM}$ is the smallest quasivariety containing both $\alg{Z}^*$ and $\alg{Z}$.

For each positive integer $n$, let $\alg{S}_{2n}$ denote the subalgebra of $\alg{Z}^*$ with universe $\{-n,\dots,-1,1,\dots,n\}$ and,
for $n\in\omega$, let $\alg{S}_{2n+1}$ be the subalgebra of $\alg{Z}$ with universe $\{-n,\dots,-1,0,1,\dots,n\}$.
The results cited above yield:
\begin{theorem}\label{thm:FG-SI-SM}
Up to isomorphism,
the algebras\/ $\alg{S}_n$ \textup{($1<n\in\omega$)} are precisely the finitely generated SI\, Sugihara monoids, whence the algebras\/ $\alg{S}_{2n+1}$ \textup{(${0<n\in\omega}$)}
are just the finitely generated SI odd Sugihara monoids.

Consequently, for each $m \in \omega$, a totally ordered $m$-generated Sugihara monoid has at most $2m +2$ elements. The bound reduces to $2m+1$ in the odd case.
\end{theorem}

%\smallskip

An element $a$ of an [I]RL $\alg{A}$ will be called \emph{negative}
if $a \leqslant e$. 
We define
\[A^{-} \seteq \{a \in A : a \leqslant e\}.\]
We say that an [I]RL $\alg{A}$ is \emph{negatively generated} 
when it is generated by negative elements, i.e., $A=\Sg^\alg{A} A^-$.
As surjective homomorphisms always map generating sets onto generating sets, the following
lemma applies.

\begin{lemma}\label{lem:neg-cone-gen-in-images}
If\/ $h\colon\alg{A}\rightarrow\alg{B}$ is a surjective homomorphism of [I]RLs and\/ $\alg{A}$ is negatively generated 
then so is\/ $\alg{B}$\textup{.}
\end{lemma}

The Sugihara monoid $\alg{Z}^{*}$ satisfies the equation
\begin{equation}
\label{eq:sigma-SM}
x = (x \wedge e) \bdot \neg (\neg x \wedge f),
\end{equation}
because De Morgan's laws reduce $\neg (\neg x \wedge f)$ to $x \vee e$, and every element of $\alg{Z}^{*}$ is comparable with $e$.
Since $\cl{SM} = \Vop(\alg{Z}^{*})$, \emph{every} Sugihara monoid $\alg{A}$ satisfies (\ref{eq:sigma-SM}) and is therefore negatively generated, as $a \wedge e \leqslant e$ and $\neg a \wedge f \leqslant f \leqslant e$ for all $a \in A$ (by Theorem~\ref{thm:idempotence}).

\begin{theorem}[{\cite[Thm.~8.5]{BMR17}}]
\label{thm:SM-Es}
Every variety of Sugihara monoids has surjective epimorphisms.
\end{theorem}
The same is true of all varieties of \emph{positive Sugihara monoids}
(i.e., RL-subreducts of Sugihara monoids) \cite[Thm.~8.6]{BMR17}.

An [I]RL is said to be
\emph{integral}
if $e$ is its greatest element.
Integral De Morgan monoids are just Boolean algebras in which $\bdot$ duplicates $\wedge$.  In the non-involutive case, integrality is less restrictive.
An integral Dunn monoid $\alg{A}$ is called a \emph{Brouwerian algebra}; it is normally
identified with its reduct $\langle A;\wedge,\vee,\to,e\rangle$, because it satisfies $x \bdot y = x \wedge y$, by (\ref{eq:meet=fusion}).

The variety of all Brouwerian algebras has the strong ES property
(see Definition~\ref{def:other-properties}(\ref{def:other-properties:strong-ES})).
The same applies to the variety of semilinear Brouwerian algebras---a.k.a.\ \emph{relative Stone algebras}.
For the origins of these results, see \cite{Kre60,EG81} and, for more comprehensive findings, Maksimova \cite{Mak03}.
More recently, it was shown in \cite[Cor.~5.7]{BMR17} that \emph{every} variety of relative Stone algebras has the (unqualified) ES property.

All varieties mentioned in and after Theorem~\ref{thm:SM-Es} consist of negatively generated Dunn/De Morgan monoids and, apart from the variety of Brouwerian algebras, their members are semilinear.
We shall show in the next two sections that
negative generation and semilinearity are enough to guarantee that a variety of Dunn/De Morgan monoids has the ES property.

\section{Semilinear Dunn Monoids}
\label{sec:no-involution}

In this section we eschew involution and consider varieties of Dunn monoids.
We start by recalling a representation theorem for totally ordered \emph{idempotent} RLs from \cite{GJM19}, and we characterize the homomorphisms between such algebras.  Our focus on the idempotent case turns out not to be restrictive, because Theorem~\ref{thm:SLDunn-neg-cone} will show that negatively generated semilinear Dunn monoids are in fact idempotent (although the same is not true for De Morgan monoids).

\smallskip
The following abbreviations are useful when working with idempotent RLs:
\[ x^{*} \seteq x \to e \;\text{ and }\; |x| \seteq x \to x. \]
In the Sugihara monoid $\alg{Z}^{*}$, the term operation $|x|$ coincides with the natural absolute value operation.
By (\ref{eq:x-y-law}), (\ref{eq:triple-equality}) and (\ref{eq:e-laws}), every [I]RL satisfies
\begin{equation}
\label{eq:star-props}
x \leqslant x^{**} \;\text{ and }\; x^{***} = x^{*} \;\text{ and }\; e \leqslant |x|.
\end{equation}
If an RL is idempotent, then it also satisfies
\begin{eqnarray*}
\label{eq:ltmod} x \leqslant |x|,\\
\label{eq:abovee} x = |x| \;\iff\; e \leqslant x,\\
\label{eq:belowe} x^{*} = |x| \;\iff\; x \leqslant e,\\
\label{eq:eqe} x = x^{*} \;\iff\; x = e.
\end{eqnarray*}

The following theorem shows that
the fusion of a totally ordered idempotent RL $\alg{A}$ resembles that of a Sugihara monoid, and that $\alg{A}$
is determined by its reduct $\langle A; \wedge,\vee,\,^{*} \rangle$, and also by its reduct $\langle A; \wedge,\vee,\left|\blank\right| \rangle$.

\begin{theorem}[{\cite[Thms.~12, 14]{Raf07}}]
\label{thm:SLidemRL-stucture}
Let $\alg{A}$ be a totally ordered idempotent RL. Then
$\alg{A}$ satisfies
\begin{equation}
\label{eq:SLidemRL-stucture}
x \bdot y = \begin{cases}
x &\text{if }\; |y| < |x|;\\
y &\text{if }\; |x| < |y|;\\
x \wedge y &\text{if }\; |x|=|y|,
\end{cases}
\;\text{ \,and\, }\;
x \to y = \begin{cases}
x^{*} \vee y &\text{if }\; x \leqslant y; \\
x^{*} \wedge y &\text{if }\; x > y.
\end{cases}
\end{equation}
\end{theorem}

Let $\alg{A}$ be a totally ordered idempotent RL.
Then
\[A^{**} \seteq \{ a^{**} : a \in A \}\]
is the universe of a subalgebra $\alg{A}^{**}$ of $\alg{A}$ which, moreover,
is termwise equivalent to
a (totally ordered) odd Sugihara monoid,
where
$\neg x \seteq x^{*}$
\cite[Lem.~3.3, Prop.~3.4]{GJM19}.
For every $c \in A^{**}$, the set
\[
A_{c} \seteq \{a \in A : a^{**} = c\}
\]
is an interval of $\alg{A}$ with greatest element $c$ \cite[Prop.~3.4]{GJM19}.
For any $\alg{A}$ as above, we define \label{pg:posets}
\[
\mathcal{A} \seteq \{\langle A_{c}\,;\, {\leqslant}|_{A_c} \rangle : c \in A^{**}\}.
\]

On the other hand, suppose $\alg{S}$ is a totally ordered odd Sugihara monoid and let
\[
\mathcal{X} = \{\langle X_c\,; \,\leqslant_c \rangle : c \in S\}
\]
be an $S$-indexed family of disjoint chains
such that each $c \in S$ is the greatest element of $X_c$.
For all $a,b \in S$ with $x \in X_a$ and $y \in X_b$, we define
\[ x \preccurlyeq y \text{ iff } a < b \text{ or }(a=b \text{ and } x \leqslant_a y).\]
Thus, $\preccurlyeq$ is the lexicographic total order on
$ S \otimes \mathcal{X} \seteq \bigcup \{X_c : c \in S\}.$
We let $\wedge$ and $\vee$ denote the meet and join operations for $\preccurlyeq$ and define 
\[ 
\alg{S} \otimes \mathcal{X} \seteq \langle S \otimes \mathcal{X}; \bdot, \to, \wedge, \vee, e \rangle, 
\]
where for $a,b \in S$ and $x \in X_a$, $y \in X_b$,
\[
x \bdot y =
\begin{cases}
x \wedge y &\text{if } a=b \leqslant e; \\
x \vee y &\text{if } e < a = b; \\
x &\text{if } a \neq b \text{ and } a \bdot^{\alg{S}} b = a; \\
y &\text{if } a \neq b \text{ and } a \bdot^{\alg{S}} b = b,
\end{cases}
\;\text{ and }\;
x \to y =
\begin{cases}
a^{*} \vee y &\text{if } x \leqslant y;\\
a^{*} \wedge y &\text{if } y < x.
\end{cases}
\]
\indent
Recall that $a \bdot b \in \{a,b\}$ for all elements $a,b$ of the Sugihara monoid $\alg{Z}^{*}$. This property is expressible as a positive universal sentence, so it holds for every totally ordered Sugihara monoid, by J\'{o}nsson's Theorem.
The above definition of $\bdot$ is therefore exhaustive.
The following representation theorem for totally ordered idempotent RLs from \cite{GJM19}
has an antecedent in \cite{Raf07}.

\begin{theorem}[{\cite[Thm.~3.5]{GJM19}}]
\label{thm:SLidemRL-representation}
For $\alg{S}$ and $\mathcal{X}$ as above,
the algebra $\alg{S} \otimes \mathcal{X}$ is a totally ordered idempotent RL satisfying $\alg{S} = (\alg{S} \otimes \mathcal{X})^{**}$ and $(S \otimes \mathcal{X})_c = X_c$ for every $c \in S$\textup. Moreover, every totally ordered idempotent RL $\alg{A}$ has this form, because
$\alg{A} = \alg{A}^{**} \otimes \mathcal{A}$.
\end{theorem}

The next two lemmas will assist in proving a characterization of homomorphisms between totally ordered idempotent RLs (Theorem~\ref{thm:SL-homs}).

\begin{lemma}[{\cite[Prop.~2.5]{Ols12}}]
\label{lem:Olson-onto-e}
Let $\alg{A}$ be a totally ordered idempotent RL and let $F \in \DFil\,\alg{A}$. If $a$ and $b$ are \emph{distinct} elements of $A$ such that $a/F = b/F$, then $b/F = e/F$.
\end{lemma}

\begin{proof}
By definition, $a/F = b/F$ means that $a \to b,\; b\to a \in F$.
By
symmetry, we may assume that $a < b$.
Then $e \not\leqslant b \to a$, by (\ref{eq:e-order}), so $b \to a < e$, because $\alg{A}$ is totally ordered.
It follows from (\ref{eq:isotone}) that $b \to (b \to a) \leqslant b \to e$. Now,
$ b \to (b \to a) = (b \bdot b) \to a = b \to a, $
by (\ref{eq:permutation}) and the idempotence of $b$.
Therefore $b \to a \leqslant b \to e$, whence $b \to e \in F$, because $b \to a \in F$.
On the other hand, $b \not\leqslant b \to a$, because otherwise $b = b \bdot b \leqslant a$. So, $b \to a < b$.
As $b\to a \in F$, we have $e \to b = b \in F$, so $b/F = e/F$.
\end{proof}
Recall that, in an algebra with a lattice reduct, any congruence class is an interval.
Specifically, if $F$ is a deductive filter of
an [I]RL
$\alg{A}$, then the set $e/F = \{a \in A : e \to a,\, a \to e \in F\} = \{a \in A : a,a^{*} \in F\}$ is an interval subuniverse of $\alg{A}$ (see \cite{HRT02} or \cite[Thm.~4.47]{GJKO07}).
When $\alg{A}$ is totally ordered, then $e/F$ is the convex closure of $\{a : e \geqslant a \in F\} \cup \{a^{*} : e \geqslant a \in F\}$, because if $e < a \in A$, then $a \to e < e$ and $a \leqslant (a \to e)^{*}$.

\begin{lemma}
\label{lem:intervals}
Let $\alg{A}$ be a totally ordered idempotent RL and
let $I$ be an interval of\/ $\alg{A}$, containing $e$, that is closed under $^{*}$\textup. Define
\[I_{*} \seteq \{a \in A : a \notin I \text{ and } a^{*} \in I\}.\]
Then
\begin{enumerate}
\item\label{lem:intervals:subalgebra} $I$ is a subuniverse of\/ $\alg{A}$\textup;
\item\label{lem:intervals:not-A**} $I_{*} \cap A^{**} = \emptyset$\textup;
\item\label{lem:intervals:below} every element of\/ $I_{*}$ is strictly below every element of\/ $I$\textup;
\item\label{lem:intervals:non-empty} if\/ $b \in I_{*}$ then $b^{*}$ is the greatest element of\/ $I$\textup;
\item\label{lem:intervals:interval} $I \cup I_{*}$ is an interval of\/ $\alg{A}$ that is closed under $^{*}$\textup.
\end{enumerate}
\end{lemma}
\begin{proof}
Item (\ref{lem:intervals:subalgebra}) holds because
$\{a \bdot b, a \to b\} \subseteq \{a,b,a^{*},b^{*}\}$ for any $a,b \in I$, by Theorem~\ref{thm:SLidemRL-stucture} (and since $e \in I$, by assumption).
Item (\ref{lem:intervals:not-A**}) holds because, otherwise, $a^{**} \in I_{*}$ for some $a \in A$, but then 
$a^{**} = a^{****} \in I$, a contradiction.

Let $b \in I_{*}$. Then $b^{*} \in I$ and so $b^{**} \in I$.
Suppose, with a view to contradiction, that $a \leqslant b$ for some $a \in I$. Then $a \leqslant b \leqslant b^{**}$, by (\ref{eq:star-props}), so because $I$ is an interval,
$b \in I$, a contradiction. Therefore, (\ref{lem:intervals:below}) holds.

For (\ref{lem:intervals:non-empty}), suppose $a > b^{*}$ for some $a \in I$.
If $b \leqslant a^{*}$, then $a \leqslant a^{**} \leqslant b^{*}$, contrary to the supposition, so $a^{*} < b \leqslant b^{**}$.
Then $b \in I$,
a contradiction.

To show (\ref{lem:intervals:interval}), notice that $I \cup I_{*}$ is clearly closed under $^{*}$, so it remains to show that $I \cup I_{*}$ is an interval. If $I_{*} = \emptyset$ we are done, so let $b$ be an arbitrary element of $I_{*}$. For any $a \in I_{*}$, we have $a^{*} = b^{*}$, by (\ref{lem:intervals:non-empty}), so $a \in A_{b^{**}}$. It follows that $I \cup I_{*}$ is the union of the overlapping intervals $I$ and $A_{b^{**}}$.
\end{proof}

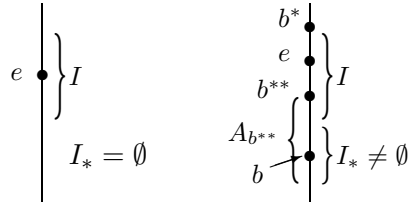
\begin{center}
\begin{picture}(72,80)

\put(30,0){\line(0,1){50}}
\put(30,48){\circle*{4}}
\put(18,46){\small $e$}
\put(40,14){$I_{*} = \emptyset$}
\put(34,39){\scalebox{1}[3]{$\}$}}
\put(40,44){\small $I$}

\put(30,50){\line(0,1){25}}

\end{picture}
\hspace{2em}
\begin{picture}(72,80)

\put(10,39){\small $b^{**}$}
\put(30,40){\circle*{4}}
\put(-1,23){\small $A_{b^{**}}$}
\put(21,15){\scalebox{1}[3]{$\{$}}
\put(30,0){\line(0,1){50}}
\put(30,53){\circle*{4}}
\put(18,53){\small $e$}
\put(34,12){\scalebox{1}[2]{$\}$}}
\put(40,14){\small $I_{*} \neq \emptyset$}
\put(34,39){\scalebox{1}[3]{$\}$}}
\put(40,44){\small $I$}

\put(30,50){\line(0,1){25}}
\put(30,66){\circle*{4}}
\put(18,66){\small $b^{*}$}

\put(30,17){\circle*{4}}
\put(15,13){\vector(3,1){12}}
\put(8,7){\small $b$}

\end{picture}
\end{center}

\begin{theorem}
\label{thm:SL-homs}
Let $\alg{A}$ and $\alg{B}$ be totally ordered idempotent RLs. A function $h \colon A \rightarrow B$ is a homomorphism from $\alg{A}$ to $\alg{B}$ iff the following hold:
\begin{enumerate}
\item\label{thm:SL-homs:interval} The set $I = h^{-1}[\{e\}]$ is an interval of $\alg{A}$\textup, which contains $e$ and is closed under $^{*}$\textup.
\item\label{thm:SL-homs:star} $h$ is an order embedding of $I_{*}$ into $B_{e} \setminus \{e\}$\textup.
\item\label{thm:SL-homs:Sugihara} $h$ is an order embedding of $A^{**} \setminus I$ into $B^{**} \setminus \{e\}$, preserving $^{*}$\textup.
\item\label{thm:SL-homs:partition} For every $a \in A^{**} \setminus I$\textup, $h$ is an order embedding of $A_a$ into $B_{h(a)}$\textup.
\end{enumerate}
\end{theorem}

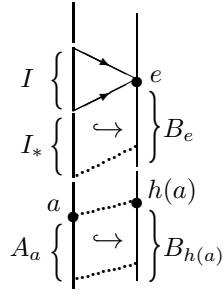
\begin{center}

\begin{picture}(90,105)(0,-5)

\put(30,85){\line(0,1){15}}
\put(10,68){\small $I$}
\put(21,65){\scalebox{1}[2]{$\{$}}
\put(30,59){\line(0,1){24}}
\put(30,83){\line(2,-1){24}}
\put(30,83){\vector(2,-1){14}}
\put(30,59){\line(2,1){24}}
\put(30,59){\vector(2,1){14}}
\put(10,43){\small $I_{*}$}
\put(37,49){$\hookrightarrow$}
\put(21,40){\scalebox{1}[2]{$\{$}}
\put(30,34){\line(0,1){24}}
\multiput(30,34)(2,1){12}{.}
\put(30,19){\circle*{4}}
\put(20,21){\small $a$}
\multiput(30,19)(2,0.5){12}{.}
\put(6,4){\small $A_{a}$}
\put(21,0){\scalebox{1}[2]{$\{$}}
\put(37,7){$\hookrightarrow$}
\put(30,-8){\line(0,1){40}}
\multiput(30,-5)(2,0.5){12}{.}

\put(54,70){\line(0,1){26}}
\put(54,70){\circle*{4}}
\put(59,70){\small $e$}
\put(58,46){\scalebox{1}[2.5]{$\}$}}
\put(64,49){\small $B_e$}
\put(54,38){\line(0,1){32}}

\put(54,24){\circle*{4}}
\put(58,26){\small $h(a)$}
\put(58,1){\scalebox{1}[2.5]{$\}$}}
\put(54,-5){\line(0,1){41}}
\put(64,4){\small $B_{h(a)}$}

\end{picture}

\end{center}

\begin{proof}
Suppose $h \colon \alg{A} \rightarrow \alg{B}$ is a homomomorphism.
The kernel $\theta$ of $h$ is $\leibniz^\alg{A}F$ for some
deductive filter $F$ of $\alg{A}$, so $\alg{A}/\theta = \alg{A} / F$.
Let $I = h^{-1}[\{e\}]$. Then $I = e/F$, which we have already noted
is an interval and a subuniverse of $\alg{A}$. In particular, $I$ is closed under $^{*}$ and contains $e$.

Let $a,b \in A \setminus I$ such that $a \neq b$.
Then $h(a) \neq h(b)$, because otherwise $a/F = b/F$, which would imply that $b/F = e/F$, by Lemma~\ref{lem:Olson-onto-e}, i.e., that $h(b) = h(e) = e$, contradicting $b \notin I$. Therefore, $h$ is injective outside of $I$.

By \cite[Lem.~3.3]{GJM19}, $h$ restricts to a homomorphism $h|_{A^{**}} \colon \alg{A}^{**} \rightarrow \alg{B}^{**}$, which in particular preserves $^{*}$. Then (\ref{thm:SL-homs:Sugihara}) holds, because $h(a) \neq e$ for any $a \in A^{**} \setminus I$.

For any $a \in I_{*}$, we have $a^{*} \in I$, so $e = e^{*} = h(a^{*})^{*} = h(a)^{**}$ and $h(a) \neq e$. Therefore, $h[I_{*}] \subseteq B_{e} \setminus \{e\}$, so (\ref{thm:SL-homs:star}) holds.
For any $a \in A^{**} \setminus I$ and $x \in A_a$, we have $x^{**} = a = a^{**}$, so $h(x)^{**} = h(a)^{**}$.
Therefore, $h[A_a] \subseteq B_{h(a)}$, so (\ref{thm:SL-homs:partition}) holds.

Conversely, let $h$ be as in the theorem
and let $I = h^{-1}[\{e\}]$.
By Theorem~\ref{thm:SLidemRL-representation}, $\alg{A} = \alg{A}^{**} \otimes \mathcal{A}$ and $\alg{B} = \alg{B}^{**} \otimes \mathcal{B}$, so the families $\{A_a : a \in A^{**}\}$ and $\{B_b : b \in B^{**}\}$ are partitions of $A$ and $B$, respectively.
Note that $I \cup I_{*} = \bigcup\{A_a : a \in I \cap A^{**}\}$, because for each $c \in A$, we have $c \in I \cup I_{*}$ iff $c^{**} \in I$ (using (\ref{eq:star-props})).
So, the sets $I$, $I_{*}$, and $A_a$ ($a \in A^{**} \setminus I$) form a partition of $A$.
It follows
from properties (\ref{thm:SL-homs:interval})--(\ref{thm:SL-homs:partition})
that $h$ is injective outside of $I$, and
that $h$ preserves order (and hence the lattice operations), in view of the definitions of $\alg{A}^{**} \otimes \mathcal{A}$ and $\alg{B}^{**} \otimes \mathcal{B}$.

Let $a \in A$. If $a \in I \cup I_{*}$ then $a^{*} \in I$, so $h(a^{*}) = e = h(a)^{*}$, by (\ref{thm:SL-homs:interval}) and (\ref{thm:SL-homs:star}). If $a \in A \setminus (I \cup I_{*})$, then $a^{**} \in A^{**} \setminus I$. From (\ref{thm:SL-homs:Sugihara}) and (\ref{thm:SL-homs:partition}), it follows that $h(a^{***}) = h((a^{**})^{*}) = h(a^{**})^{*}$ and $h(a) \in h[A_{a^{**}}] \subseteq B_{h(a^{**})}$, i.e., $h(a)^{**} = h(a^{**})$. So,
$h(a^{*}) = h(a^{***}) = h(a^{**})^{*} = h(a)^{***} = h(a)^{*}$.
Therefore, $h$ preserves $^{*}$.

For $a,b \in A$, the characterization of $a\to b$ in Theorem~\ref{thm:SLidemRL-stucture} shows that preservation of $\to$ follows from that of $\wedge$, $\vee$ and~$^{*}$, except when $a > b$ but $h(a) =h(b)$.
In this situation $a,b \in I$, because $h$ is injective outside of $I$,
so $a \to b \in I$, by Lemma~\ref{lem:intervals}(\ref{lem:intervals:subalgebra}), whence 
\[ h(a) \to h(b) = e \to e = e = h(a \to b). \]
Therefore, $h$ also preserves $\left|\blank\right|$.

Again, by Theorem~\ref{thm:SLidemRL-stucture}, since $h$ preserves $\vee$, $\wedge$ and $\left|\blank\right|$, to show preservation of $\bdot$, we need only consider cases where $|a| > |b|$ but $|h(a)| = |h(b)|$ (i.e., $h(|a|)=h(|b|)$) for some $a,b \in A$. Then $|a|,|b| \in I$, so $a,b \in I \cup I_{*}$, because
$|c| \in \{c,c^{*}\}$ for all $c \in A$, by Theorem~\ref{thm:SLidemRL-stucture}.
We have $h(a \bdot b) = h(a)$, since $a \bdot b = a$, and $h(a) \bdot h(b) = h(a) \wedge h(b)$, so it remains to show that $h(a) \leqslant h(b)$.

If $b \in I_{*}$ then
$b \leqslant e$, by Lemma~\ref{lem:intervals}(\ref{lem:intervals:below}), so $|b| = b^{*} \geqslant |a|$, by  Lemma~\ref{lem:intervals}(\ref{lem:intervals:non-empty}), since $|a| \in I$.  But this contradicts $|a|>|b|$, so $b\notin I_{*}$.
If $a,b \in I$, then $h(a) = e = h(b)$. Lastly, if $a \in I_{*}$ and $b \in I$, then $a \leqslant e$, so $h(a) \leqslant h(e)= e= h(b)$.
\end{proof}

The following theorem exhibits a variety of semilinear Dunn monoids with the ES property 
which has members that are not negatively generated.
(One can easily construct totally ordered idempotent RLs that are not negatively generated, using Theorem~\ref{thm:SLidemRL-representation}; also see the examples before Lemma~\ref{lem:SLid-simple}.)

\begin{theorem}
\label{thm:SLid-ES}
Epimorphisms are surjective in the variety of all idempotent semilinear RLs.
\end{theorem}

\begin{proof}
Let $\alg{B}$ be a proper subalgebra of a totally ordered idempotent RL $\alg{A}$.
Let $a \in A \setminus B$.
We shall show that $\alg{B}$ is not epic in $\alg{A}$ by constructing a totally ordered idempotent RL $\alg{C}$ and two homomorphisms from $\alg{A}$ into $\alg{C}$ which agree on $\alg{B}$ but differ at $a$.
It then follows from Theorem~\ref{thm:Campercholi} that the variety of all idempotent semilinear RLs has the ES property.
We split into two cases: $a \in A^{**}$ and $a \notin A^{**}$.

First suppose that $a \in A^{**}$.
Without loss of generality, $a < e$. Indeed, if $e \leqslant a$, then $e \geqslant a^{*} = a^{***} \in A^{**}$; moreover, $a^{*} \notin B$, because otherwise, $a = a^{**} \in B$.

Now $F \seteq \{b \in A : b > a\}$ is a deductive filter of $\alg{A}$. Let $q$ be the canonical surjection from $\alg{A}$ to the totally ordered algebra $\alg{C} \seteq \alg{A}/F$.
We use the notation $[x,y]$ to denote the interval from $x$ to $y$, i.e., the set $\{z : x \leqslant z \leqslant y\}$.
Consider the set
\[I \seteq [a,a^{*}] \cup [a,a^{*}]_{*} = [a,a^{*}]\cup\{x \in A : x \notin [a,a^{*}] \text{ and } x^{*} \in [a,a^{*}]\}.\]
Define a map $h\colon A \rightarrow A/F$ by
\[h(x) = \begin{cases}
e^{\alg{C}} &\text{if } x \in I; \\
q(x) &\text{otherwise}.
\end{cases}\]
Note that $h(a) = e^{\alg{C}}$, since $a \in I$,  so $h(a) \neq q(a)$ (because, if $a \in e/F$ then $a = e \to a \in F$, which is not the case). We now show that $h$ is a homomorphism.
As $[a,a^{*}]$ is an interval containing $e$ that is closed under $^{*}$ (because $a \in A^{**}$), the same of true of $I$, by Lemma~\ref{lem:intervals}(\ref{lem:intervals:interval}).
Furthermore, $h^{-1}[\{e^{\alg{C}}\}] = I$, because $q^{-1}[\{e^{\alg{C}}\}] = e/F \subseteq \{b \in A : a < b \leqslant a^{*}\} \subseteq I$.
So, condition (\ref{thm:SL-homs:interval}) of Theorem~\ref{thm:SL-homs} holds.
Note that $a^{*} = \max I$, by Lemma~\ref{lem:intervals}(\ref{lem:intervals:below}).

If $b \in I_{*}$ then $b \notin I$ and $b^{*}$ is the greatest element of $I$, by Lemma~\ref{lem:intervals}(\ref{lem:intervals:non-empty}), so $b^{*} = a^{*}$. But, since $b \notin [a,a^{*}]$ and $b^{*} \in [a,a^{*}]$, we have $b \in [a,a^{*}]_{*}$, a contradiction. So, $I_{*} = \emptyset$, and condition~(\ref{thm:SL-homs:star}) of Theorem~\ref{thm:SL-homs} is vacuously satisfied.

As $q$ is a homomorphism between totally ordered idempotent RLs, Theorem~\ref{thm:SL-homs} applies to $q$. In particular, the following conditions hold:
\begin{itemize}
\item[(\ref{thm:SL-homs:Sugihara})] $q$ is a $^{*}$-preserving order embedding of $A^{**} \setminus (e/F)$ into $C^{**} \setminus \{e\}$;
\item[(\ref{thm:SL-homs:partition})] for every $a \in A^{**} \setminus (e/F)$, $q$ is an order embedding of $A_a$ into $C_{q(a)}$.
\end{itemize}
As $e/F \subseteq I$, conditions (\ref{thm:SL-homs:Sugihara}) and (\ref{thm:SL-homs:partition}) also hold for $h$.
So, $h$ is a homomorphism, by Theorem~\ref{thm:SL-homs}.

To show that $h|_{B} = q|_{B}$, we let $b \in B \cap I$ and prove that $q(b) = e^{\alg{C}}$, i.e., that $b,b^{*} \in F$. Note that $b^{*} \in [a,a^{*}]$. If $b^{*} \notin F$ then $b^{*} \leqslant a$ by the definition of $F$, so $a = b^{*} \in B$, a contradiction. Therefore $b^{*} \in F$, as claimed.
Suppose that $b \in [a,a^{*}]_{*}$. By Lemma~\ref{lem:intervals}(\ref{lem:intervals:non-empty}), $b^{*} = a^{*}$, so $a = a^{**} = b^{**} \in B$, a contradiction. So, $b \in [a,a^{*}]$. Since $a \neq b$, we get $a < b$, i.e., $b \in F$, completing the proof that $h|_{B} = q|_{B}$.

Now suppose that $a \notin A^{**}$.
Let $\alg{A}'_s = \langle A_s;\, {\leqslant}|_{A_s} \rangle$ whenever $a^{**} \neq s \in A^{**}$.
Define $\alg{A}'_{a^{**}} = \langle A'_{a^{**}};\, \leqslant' \rangle$, where $A'_{a^{**}} = A_{a^{**}} \cup \{c\}$ for some fresh element $c \notin A$ and $\leqslant'$ is the total order on $A'_{a^{**}}$ that extends ${\leqslant}|_{A_{a^{**}}}$ with $c <' a$ and $b <' c$ whenever $a > b \in A_{a^{**}}$.
Let $\mathcal{A}' = \{\alg{A}'_s : s \in A^{**}\}$ and $\alg{C} = \alg{A}^{**} \otimes \mathcal{A}'$. By Theorem~\ref{thm:SLidemRL-representation}, $\alg{C}$ is a totally ordered idempotent RL.
By Theorem~\ref{thm:SL-homs}, the inclusion map $i\colon \alg{A} \rightarrow \alg{C}$ is a homomorphism, and so is the map
\[h \colon x \mapsto \begin{cases}
c &\text{if } x = a;\\
x &\text{otherwise}.
\end{cases}.\]
Note that $h$ and $i$ differ only at $a$, so $h|_{B} = i|_{B}$.
\end{proof}

It was recently shown in \cite[Thm.~6.6]{GJM19} that the variety of semilinear idempotent RLs has the amalgamation property. Combining this observation with Theorems~\ref{thm:SLid-ES} and \ref{thm:linkBetweenESandAmalgamation}, we obtain:
\begin{corollary}
\label{cor:indemSemiCRLhasStrongESandStrongAmalgamation}
The variety of semilinear idempotent RLs has the strong amalgamation property, and hence the strong ES property.
\end{corollary}

The proof of Theorem~\ref{thm:SLid-ES} essentially showed that the class of \emph{totally ordered} idempotent RLs has the strong ES property.
Nevertheless, we could not have deduced from
this
alone that the whole variety of \emph{semilinear} idempotent RLs
has the strong ES property, because
Theorem~\ref{thm:Campercholi} has no analogue for the strong ES property.

To see this, let $\cl{L}^4$ denote the variety generated by the four-element totally ordered Brouwerian
algebra.
The strong ES property holds for $\cl{L}^4_\textup{FSI}$, but fails for $\cl{L}^4$.
Indeed, Maksimova showed in \cite[Thm.~4.3]{Mak03} that
just six nontrivial varieties of Brouwerian algebras have
the strong ES property, and only three of these consist of semilinear algebras, namely the class of all relative Stone algebras and the varieties generated, respectively, by the two-element and the three-element relative Stone algebras.
That $\cl{L}^4_\textup{FSI}$ has the strong ES property can be deduced from the proof of Theorem~\ref{thm:GSM-ES} below.

Not all varieties of semilinear idempotent RLs have the ES property, as we shall see in Example~\ref{ex:SLid-ES-failure}. But 
we shall prove in Theorem~\ref{thm:GSM-ES}
that epimorphisms are surjective in all varieties of \emph{negatively generated} semilinear idempotent RLs. 

\begin{definition}
\label{def:GSM}
The variety $\cl{GSM}$
of \emph{generalized Sugihara monoids}
consists of the semilinear idempotent RLs that satisfy
\begin{equation*}
\label{eq:GSM}
(x \vee e)^{**} =  x \vee e,
\end{equation*}
or equivalently, $e \leqslant x \implies x^{**} =  x$.
\end{definition}
The main significance of $\cl{GSM}$ lies in the next theorem. 

\begin{theorem}[{\cite[Cor.~3.5]{GR15}}]
\label{thm:GSM-neg-cone}
A semilinear idempotent RL is a generalized Sugihara monoid iff it is negatively generated. \end{theorem}
\noindent
In the proof of this theorem, one uses the fact that all generalized Sugihara monoids satisfy
\begin{equation}
\label{eq:sigma}
x =  (x \wedge e) \bdot (x^{*} \wedge e)^{*}.
\end{equation}

\begin{corollary}[{\cite{GR15}}] 
\label{cor:GSM-representation}
A totally ordered idempotent RL $\alg{A}$ is a generalized Sugihara monoid iff\/ $A_{c} = \{c\}$ for every $c > e$\textup{.}
\end{corollary}
\begin{proof}
($\Rightarrow$):  Let $e < c \in A$.  As $\alg{A} \in \cl{GSM}$, we have $c^{**} = c$, 
so $c \in A_c$.  Now let $d \in A_c$, i.e., $d^{**} = c$, so $d^{*} = d^{***} = c^{*} \leqslant e$.  We must show that $d = c$.  If $d \leqslant e$, then $e \leqslant d^{*}$, so $d^{*} = e$, whence $c = d^{**} = e$, a contradiction.  Consequently, $e < d$.  Then, since $\alg{A} \in \cl{GSM}$, we have $d = d^{**} = c$.

($\Leftarrow$):  Suppose $A_c = \{c\}$ whenever $e < c \in A$.  To see that $\alg{A} \in \cl{GSM}$, let $e \leqslant a \in A$.  If $a = e$, then $a^{**} = a$, so suppose $e < a$.  Then $e < a^{**}$ (as $a \leqslant a^{**}$), so $A_{a^{**}} = \{a^{**}\}$, by assumption.  But $a \in A_{a^{**}}$, so $a = a^{**}$, as required.
\end{proof}

The next theorem strengthens \cite[Thm.~13.1]{GR15},
which stated that every variety of generalized Sugihara monoids has the weak ES property.
It also unifies two findings from \cite{BMR17}: all varieties of positive Sugihara monoids and all varieties of relative Stone algebras have the ES property.

\begin{theorem}
\label{thm:GSM-ES}
All varieties of generalized Sugihara monoids have surjective epimorphisms.
\end{theorem}

\begin{proof}
Assume, with a view to contradiction, that $\cl{K}$ is a subvariety of $\cl{GSM}$ without the ES property. Then, by Theorem~\ref{thm:Campercholi},
there exists $\alg{A} \in \cl{K}_\textup{FSI}$ (i.e., a totally ordered $\alg{A} \in \cl{K}$) with a $\cl{K}$-epic proper subalgebra $\alg{B}$.

Since $\alg{A}$ is negatively generated,
there exists $a \in A^{-} \setminus B$, so $a < e$.
Then $F \seteq \{b \in A :  a < b\}$ is a deductive filter of $\alg{A}$. Let $\alg{C} \seteq \alg{A}/F$, and let $q \colon \alg{A} \rightarrow \alg{C}$ be the canonical surjection. Note that $\alg{C}$ is totally ordered and $\alg{C} \in \cl{K}$, because $\cl{K}$ is a variety.

Recall that $a \leqslant a^{**}$. 
If $a = a^{**}$, then $a \in A^{**}$ and we can use the 
first homomorphism in the proof of
Theorem~\ref{thm:SLid-ES} to show that $\alg{B}$ is not $\cl{K}$-epic in $\alg{A}$, a contradiction.

So, we may suppose that $a < a^{**}$. In this case, define $h \colon A \rightarrow C$ by
\[h(x) = \begin{cases}
e^{\alg{C}} &\text{ if } x = a; \\
q(x) &\text{ otherwise}.
\end{cases}
\]
Then $h^{-1}[\{e^{\alg{C}}\}] = (e/F) \cup \{a\}$. We claim that $(e/F) \cup \{a\} = [a,a^{*}]$, which is clearly an interval of $\alg{A}$ containing $e$ that is closed under $^{*}$.
If $b \in [a,a^{*}]$ and $b \neq a$, we must show that $b/F = e/F$, i.e., that $a < b,b^{*}$. Clearly $a <b$ and $a \leqslant b^{*}$. If $a = b^{*}$, then $a^{**} = b^{***} = b^{*} =a$, contradicting the assumption that $a < a^{**}$. So, $a < b^{*}$, as required.

Because $q$ satisfies conditions (\ref{thm:SL-homs:star})--(\ref{thm:SL-homs:partition}) of Theorem~\ref{thm:SL-homs}, and $q^{-1}[\{e^{\alg{C}}\}] = e/F \subseteq h^{-1}[\{e^{\alg{C}}\}]$, it is easy to see that $h$ satisfies the conditions of Theorem~\ref{thm:SL-homs}.
So, $h$ is a homomorphism from $\alg{A}$ to $\alg{C}$.
Clearly, $h|_{B} = q|_{B}$, but $h(a) \neq q(a)$. Therefore, $\alg{B}$ is not $\cl{K}$-epic in $\alg{A}$, a contradiction.
\end{proof}

Recall that a filter $F$ of a lattice $\langle L;\wedge,\vee\rangle$
is said to be \emph{prime}
if its complement $L\setminus F$ is closed under the binary operation $\vee$.
We say that a square-increasing [I]RL $\alg{A}$ has \emph{infinite depth} if its poset of prime deductive filters contains an infinite descending chain; otherwise it has \emph{finite depth}. This definition is equivalent to the one employed in \cite{MRW4}, where a slightly stronger version of the following result is proved.

\begin{theorem}[{\cite[Thm.~8.1]{MRW4}}]
\label{thm:ES-finite-depth}
Let\/ $\cl{K}$ be a variety of square-increasing [I]RLs,
such that each FSI member of\/ $\cl{K}$
is
negatively generated and
has finite depth.
Then every\/ $\cl{K}$-epimorphism is surjective.
\end{theorem}

Every variety $\cl{K}$ of RLs with the ES property exhibited thus far in this paper has at least one of the following two properties: (i) $\cl{K}$ is generated by algebras that are negatively generated (as in
Theorems~\ref{thm:ES-finite-depth} and
\ref{thm:GSM-ES}), or (ii) $\cl{K}$ has infinite depth (as in Theorems~\ref{thm:SLid-ES} and \ref{thm:GSM-ES}).
In Theorem~\ref{thm:SLid-simple-ES} below,
we identify a variety of semilinear Dunn monoids with surjective epimorphisms which satisfies neither (i) nor (ii).

Let $\alg{2}^{+}$ denote the two-element Brouwerian algebra. Recall that the three-element Sugihara monoid $\alg{S}_3$ has universe $\{-1,0,1\}$.
For any chain $\alg{P}$ with greatest element $1$, we abbreviate $\alg{S}_3 \otimes \{\{-1\},\{0\},\alg{P}\}$ as $\alg{S}_3 \oplus \alg{P}$.

\begin{center}
\begin{tabular}{ccc}

\begin{picture}(35,45)(-10,0)

\put(12,36){\circle*{4}}
\put(12,36){\line(0,-1){10}}

\put(12,24){\line(0,-1){4}}
\put(12,18){\line(0,-1){4}}
\put(12,12){\line(0,-1){4}}

\put(17,32){\small $1$}
\put(-9,17){\small $\alg{P}$:}
\end{picture}
&
\begin{picture}(45,40)(-14,0)

\put(12,27){\circle*{4}}
\put(12,27){\line(0,-1){9}}
\put(12,18){\circle*{4}}
\put(12,18){\line(0,-1){9}}
\put(12,9){\circle*{4}}

\put(20,26){\small $1$}
\put(20,14){\small $0$}
\put(16,3){\small $-1$}
\put(-13,17){\small $\alg{S}_3$:}
\end{picture}
&
\begin{picture}(66,50)(-33,-10)

\put(12,36){\circle*{4}}
\put(12,36){\line(0,-1){10}}

\put(12,24){\line(0,-1){4}}
\put(12,18){\line(0,-1){4}}
\put(12,12){\line(0,-1){4}}

\put(12,2){\circle*{4}}
\put(12,2){\line(0,-1){9}}
\put(12,-6){\circle*{4}}

\put(20,30){\small $1$}
\put(20,-0){\small $0$}
\put(16,-10){\small $-1$}
\put(-35,7){\small $\alg{S}_3 \oplus \alg{P}$:}
\end{picture}

\end{tabular}
\end{center}

\begin{lemma}[{\cite[Thm.~3.7]{Ols08}}]
\label{lem:SLid-simple}
A semilinear idempotent RL is simple iff it is isomorphic to
$\alg{2}^+$ or $\alg{S}_3 \oplus \alg{P}$ for some chain $\alg{P}$ with top element\/ $1$\textup.
\end{lemma}

Let $\cl{S}$ be the class of all simple totally ordered idempotent RLs.
\begin{theorem}
\label{thm:SLid-simple-ES}
Epimorphisms are surjective in $\Vop(\cl{S})$.
\end{theorem}

\begin{proof}
Let $\alg{A} \in \Vop(\cl{S})_\textup{FSI}$ and let $\alg{B}$ be a proper subalgebra of $\alg{A}$,
so $\alg{A}$ is nontrivial.  Just as in Theorems~\ref{thm:SLid-ES} and~\ref{thm:GSM-ES}, we must show that $\alg{B}$ is not
$\Vop(\cl{S})$-epic in $\alg{A}$.

By J\'{o}nsson's Theorem, the FSI members of $\Vop(\cl{S})$ belong to
$\Hop\Sop\Puop(\cl{S})$, but the criterion for simplicity in Lemma~\ref{lem:RL-properties}(\ref{lem:RL-properties:simple})
is first order-definable and therefore persists in ultraproducts (by \L os' Theorem \cite[Thm.~5.21]{Ber12}), 
while the CEP ensures that nontrivial subalgebras of simple algebras are simple. 
Therefore, $\alg{A}$ is simple, since $\alg{A} \in \Vop(\cl{S})_\textup{FSI}$.

Thus, $\alg{A}$ is isomorphic to $\alg{2}^{+}$ or $\alg{S}_3 \oplus \alg{P}$ for some chain $\alg{P}$ with greatest element $1$,
by Lemma~\ref{lem:SLid-simple}.

In the first case, $B = \{e\}$. The identity map from $\alg{A}$ to itself, and the map sending $\alg{A}$ onto the trivial subalgebra of $\alg{A}$, are different homomorphisms that agree on $B$. So, $\alg{B}$ is not $\Vop(\cl{S})$-epic in $\alg{A}$.

We may therefore suppose that $\alg{A} = \alg{S}_3 \oplus \alg{P}$.
If $\alg{B}$ is trivial, we are done, as in the previous paragraph. So, we may assume that $\alg{B}$ is nontrivial, in which case $S_3 \subseteq B$.
Let $c \in A \setminus B$. Then $c \in P \setminus \{1\}$.

As in Theorem~\ref{thm:SLid-ES}, let $P'=P\cup\{d\}$ for some fresh element $d \notin A$ and extend the total order of $\alg{P}$ to $\alg{P}'$ by defining $d$ to be the immediate predecessor of $c$. By Theorem~\ref{thm:SL-homs}, the inclusion map $i$ from $\alg{S}_3 \oplus \alg{P}$ to $\alg{S}_3 \oplus \alg{P}' \in \cl{S}$ is a homomorphism, and the map
\[h \colon x \mapsto \begin{cases}
d &\text{if }x = c;\\
x &\text{otherwise}
\end{cases}\]
is also a homomorphism, differing from $i$ only at $c$, so that $h|_{B} = i|_{B}$.
Thus, $\alg{B}$ is not $\Vop(\cl{S})$-epic in $\alg{A}$.
\end{proof}

We now exhibit a subvariety of $\Vop(\cl{S})$ which does not have the ES property.
Let $\mathit{2}$ be the two-element chain with elements $c < 1$.
\begin{center}
\begin{picture}(64,40)(-34,0)

\put(12,36){\circle*{4}}
\put(12,36){\line(0,-1){9}}
\put(12,27){\circle*{4}}
\put(12,27){\line(0,-1){9}}
\put(12,18){\circle*{4}}
\put(12,18){\line(0,-1){9}}
\put(12,9){\circle*{4}}

\put(20,33){\small $1$}
\put(20,24){\small $c$}
\put(20,12){\small $0$}
\put(16,2){\small $-1$}
\put(-35,20){\small $\alg{S}_3 \oplus \mathit{2}$:}
\end{picture}
\end{center}

\begin{example}
\label{ex:SLid-ES-failure}
$\Vop(\alg{S}_3 \oplus \mathit{2})$ does not have the ES property.
\end{example}

\begin{proof}
Let $\cl{K}=\Vop(\alg{S}_3 \oplus \mathit{2})$. We show that $\alg{S}_3$ is a $\cl{K}$-epic subalgebra of $\alg{S}_3 \oplus \mathit{2}$.

Let $g,h \colon \alg{S}_3 \oplus \mathit{2} \rightarrow \alg{C}$ be two homomorphisms into $\alg{C} \in \cl{K}_\text{SI}$
such that $g|_{S_3} = h|_{S_3}$. By Lemma~\ref{lem:FSI-image}, it suffices to show that $g = h$.

Since $\alg{S}_3 \oplus \mathit{2}$ is simple, and $g$ and $h$ agree on a non-neutral element, $g$ and $h$ are both embeddings or they both have range $\{e\}$. In the second case, clearly $g=h$. So, we assume that $g$ and $h$ are embeddings.

By J\'{o}nsson's theorem, $\alg{C} \in \Hop\Sop\Puop(\alg{S}_3 \oplus \mathit{2})$.
Since $\alg{S}_3 \oplus \mathit{2}$ is finite and simple, $\alg{C}$ is isomorphic to $\alg{S}_3$ or to $\alg{S}_3 \oplus \mathit{2}$.
Since $g$ and $h$ are embeddings, the first case is ruled out on cardinality grounds, so $\alg{C} \cong \alg{S}_3 \oplus \mathit{2}$. But then $g = h$ because $\alg{S}_3 \oplus \mathit{2}$ has no nontrivial automorphism.
\end{proof}
Note that $\alg{S}_3 \oplus \mathit{2}$ is not negatively generated (as the subuniverse generated by $\{-1,0\}$ excludes $c$).
Also, as $\alg{S}_3 \oplus \mathit{2}$ is finite and has a $\Vop(\alg{S}_3 \oplus \mathit{2})$-epic proper subalgebra, $\Vop(\alg{S}_3 \oplus \mathit{2})$ fails to have even the weak ES property.
(It follows from \cite[Cor.~6.5]{Cam18} that, in finitely generated varieties of lattice-based algebras,
the weak ES property entails the ES property, but this becomes false if we replace `finitely generated' by `locally finite' \cite[Sec.~6]{BMR17}.)

%\smallskip

We now relax the condition of idempotence and consider varieties of semilinear RLs that are merely square-increasing.

\begin{theorem}
\label{thm:idempotent-generators}
Let $\alg{A}$ be a totally ordered Dunn monoid that is generated by a set $X$ of idempotent elements. Then $\alg{A}$ is idempotent.
\end{theorem}

\begin{proof}
Let $a \in A$.
Then $a = t^{\alg{A}}(a_1, \dots, a_n)$ for some
RL-term $t(x_1, \dots, x_n)$ and some $a_1, \dots, a_n \in X$.  Let $\vec{a}$ abbreviate $a_1, \dots, a_n$.
We show that $a = a^{2}$ by induction on the complexity $\# t$ of $t$.  For brevity, we assume below that all terms are evaluated in $\alg{A}$.

When $\# t = 0$, clearly $t(\vec{a})^2 = t(\vec{a})$, because $t \in \{e, x_1, \dots, x_n\}$.

Assume that $s$ and $r$ are RL-terms with $\# s, \# r < \# t$, where $s(\vec{a})^2 = s(\vec{a})$ and $r(\vec{a})^2 = r(\vec{a})$.

If $t$ is $s \wedge r$ or $s \vee r$, then $t(\vec{a}) \in \{s(\vec{a}),r(\vec{a})\}$, since $\alg{A}$ is totally ordered, and we are done.
If $t = s \bdot r$, then $t(\vec{a})^2 = s(\vec{a})^2 \bdot r(\vec{a})^2 = s(\vec{a}) \bdot r(\vec{a}) = t(\vec{a})$, by the induction hypothesis.
Lastly, suppose that $t = s \to r$. Note that $t(\vec{a}) \leqslant t(\vec{a})^2$, since $\alg{A}$ is square-increasing.
On the other hand, by (\ref{eq:x-y-law}),
\[ s(\vec{a}) \bdot (s(\vec{a}) \to r(\vec{a}))^2 = s(\vec{a})^2 \bdot (s(\vec{a}) \to r(\vec{a}))^2 \leqslant r(\vec{a})^2 = r(\vec{a}),\]
so 
$t(\vec{a})^2 = (s(\vec{a}) \to r(\vec{a}))^2 \leqslant s(\vec{a}) \to r(\vec{a}) = t(\vec{a})$, by the law of residuation (\ref{eq:residuation}).
\end{proof}

\begin{theorem}
\label{thm:SLDunn-neg-cone}
Let $\alg{A}$ be an semilinear Dunn monoid. Then the following are equivalent:
\begin{enumerate}
\item \label{thm:SLDunn-neg-cone:neg-cone} $\alg{A}$ is negatively generated;
\item \label{thm:SLDunn-neg-cone:FSI} $\alg{A}$ is a generalized Sugihara monoid;
\item \label{thm:SLDunn-neg-cone:equation} $\alg{A}$ satisfies equation~\textup{(\ref{eq:sigma}).} 
\end{enumerate}
\end{theorem}

\begin{proof}
(\ref{thm:SLDunn-neg-cone:neg-cone}) $\Rightarrow$ (\ref{thm:SLDunn-neg-cone:FSI}):  By the
Subdirect Decomposition Theorem, $\alg{A}$ embeds into $\prod_{i \in I} \alg{A}_i$ for some set $\{\alg{A}_i : i \in I\}$ of totally ordered Dunn monoids, where each $\alg{A}_i$ is a homomorphic image of $\alg{A}$.
For each $i \in I$, we have $\alg{A}_i = \SgA^{\alg{A}_i} A_i^{-}$, by Lemma~\ref{lem:neg-cone-gen-in-images}.
By (\ref{eq:meet=fusion}), every element of $A_i^{-}$ is idempotent,
so $\alg{A}_i$ is idempotent, by Theorem~\ref{thm:idempotent-generators}. Therefore, each $\alg{A}_i \in \cl{GSM}$, by Theorem~\ref{thm:GSM-neg-cone}, so $\alg{A}$ is a generalized Sugihara monoid.

For (\ref{thm:SLDunn-neg-cone:FSI}) $\Rightarrow$ (\ref{thm:SLDunn-neg-cone:equation}), see the proof of \cite[Cor.~3.5]{GR15}.
And (\ref{thm:SLDunn-neg-cone:equation}) $\Rightarrow$ (\ref{thm:SLDunn-neg-cone:neg-cone}) follows from the form of (\ref{eq:sigma}), as
$a \wedge e$ and $a^{*} \wedge e$ belong to $A^-$, for all $a \in A$.
\end{proof}

\begin{corollary}
\label{cor:SLDun-neg-cone-loc-fin-variety}
The negatively generated semilinear Dunn monoids form a locally finite variety, namely the variety of generalized Sugihara monoids.
\end{corollary}

Indeed, it is shown in \cite[Thm.~18]{Raf07}
that the variety of semilinear idempotent RLs is locally finite, therefore $\cl{GSM}$ is as well.

The following characterization of locally finite varieties is often useful. (A proof can be found in \cite[Thm.~1]{Raf07}, for instance.)
\begin{lemma}
\label{fct:FG-function}
A variety\/ $\cl{K}$ of finite type is locally finite iff
there is a function $p\colon\omega\rightarrow\omega$ such that,
for each $n\in\omega$, every $n$-generated member of\/ $\cl{K}_\textup{SI}$ has at most $p(n)$ elements.
\end{lemma}

For each $n \in \omega$, an $n$-generated totally ordered idempotent RL has at most $3n+1$ elements \cite[Thm.~17]{Raf07}.
The bound reduces to $n+1$ in the integral case, i.e., in the subvariety of relative Stone algebras.

The following is therefore a paraphrase of Theorem~\ref{thm:GSM-ES}.

\begin{corollary}
\label{cor:SLDunn-ES}
Let\/ $\cl{D}$ be a variety of negatively generated semilinear Dunn monoids.
Then $\cl{D}$ has surjective epimorphisms.
\end{corollary}

The variety of \emph{all} semilinear Dunn monoids does not have the ES property, however. This is illustrated by the next theorem and the examples discussed after it.

\begin{theorem}
\label{thm:Aps}
Let\/ $\cl{K}$ be a variety of semilinear Dunn or De Morgan monoids containing a totally ordered algebra $\alg{A}$ which is generated by some
$a \in A$ that satisfies $a = a^{n} \to a^{n+1}$ for some positive integer $n$, where $a^{n+1}$ generates a proper subalgebra of $\alg{A}$.
Then\/ $\cl{K}$ lacks the weak ES property.
\end{theorem}

\begin{proof}
It suffices to show that $\alg{B} = \SgA^{\alg{A}} \{a^{n+1}\}$ is a
$\cl{K}$-epic subalgebra of $\alg{A}$ (see Definition~\ref{def:other-properties}(\ref{def:other-properties:weak-ES})).  Suppose, on the contrary, that
$h,g \colon \alg{A} \rightarrow \alg{C}$ are different homomorphisms that agree at $a^{n+1}$, where $\alg{C}\in\cl{K}$.
By Lemma~\ref{lem:FSI-image}, we may assume that $\alg{C} \in \cl{K}_\textup{SI}$,
so $\alg{C}$ is totally ordered.
Now $h(a) \neq g(a)$, because $\alg{A}$ is generated by $a$.  By symmetry, we may assume that ${h(a) < g(a)}$, so
 $h(a^{n}) = h(a)^{n} \leqslant g(a)^{n} = g(a^{n})$.
Then
\[g(a) \bdot h(a^{n}) \leqslant g(a) \bdot g(a^{n}) = g(a^{n+1}) = h(a^{n+1}),\]
whence $g(a) \leqslant h(a^{n}) \to h(a^{n+1}) = h(a^{n} \to a^{n+1}) = h(a)$, by the law of residuation (\ref{eq:residuation}), a contradiction.
\end{proof}

For each positive integer $p$, consider the totally ordered De Morgan monoid $\alg{A}_p$ on the chain $0 < 1 < 2 < \dots < 2^{p+1}$, where fusion is multiplication, truncated at $2^{p+1}$. For each $p>2$, the algebra $\alg{A}_p$ is generated by $2$, and $2 = 2^{p-1} \to 2^{p}$.  Also, $f = 2^{p}$, and the subalgebra $\SgA^{\alg{A}_p}\{f\}$ has universe $\{0,1,2^p,2^{p+1}\}$,
so $\alg{A}_p$ satisfies the conditions of Theorem~\ref{thm:Aps} with $n = p-1$.
When $p$ is prime, then $\alg{A}_p$ has no proper subalgebra other than $\SgA^{\alg{A}_p}\{f\}$ \cite[Ex.~9.1]{MRW2}.

An analogous situation holds for the involution-less reducts of these algebras.
For each positive integer $p$, let $\alg{A}_p^+$ denote the Dunn monoid reduct of $\alg{A}_p$.  Then $2$ still generates $\alg{A}_p^+$, and $2 = 2^{p} \to 2^{p+1}$.  As $2^{p+1}$ is idempotent in $\alg{A}_p^+$, it generates an idempotent subalgebra of $\alg{A}_p^+$, by Theorem~\ref{thm:idempotent-generators}, which must therefore be a proper subalgebra. In fact, $\Sg^{\alg{A}_p^+} \{2^{p+1}\} = \{0,1,2^{p+1}\}$. Thus, $\alg{A}_p^+$ satisfies the conditions of Theorem~\ref{thm:Aps}, with $a = 2$ and $n = p$.
When $p$ is prime, the only nontrivial proper subalgebra of $\alg{A}_{p-1}^+$ has universe $\{0,1,2^p\}$; it is isomorphic to the Dunn monoid reduct of $\alg{S}_3$.

By Lemma~\ref{lem:RL-properties}(\ref{lem:RL-properties:simple}), $\alg{A}_p$ and $\alg{A}_p^+$ are simple.
For distinct primes $p,q$, J\'{o}nsson's Theorem shows that
$\Vop(\alg{A}_{p}) \neq \Vop(\alg{A}_{q})$ and $\Vop(\alg{A}_{p-1}^+) \neq \Vop(\alg{A}_{q-1}^+)$.

Somewhat more can be said, because for any subset $\cl{X}$ of $\{\alg{A}_p:p\text{ prime}\}$ or of $\{\alg{A}_{p-1}^+:p\text{ prime}\}$, the variety $\Vop(\cl{X})$ still satisfies the conditions of Theorem~\ref{thm:Aps} and therefore lacks the weak ES property.
Moreover, De Morgan/Dunn monoids have equationally definable principal congruences (EDPC) \cite[Thm.~3.55]{GJKO07}.  In any variety $\cl{K}$ of finite type with EDPC, and for every finite algebra $\alg{A}\in\cl{K}_\text{SI}$, the class
$\cl{K}_\alg{A}\seteq\{\alg{B}\in\cl{K}:\alg{A}\notin\Sop\Hop(\alg{B})\}$
is a subvariety of $\cl{K}$, and for every subvariety $\cl{W}$ of $\cl{K}$, we have $\alg{A}\in\cl{W}$ or $\cl{W}\subseteq\cl{K}_\alg{A}$, and not both \cite{BP82},\,\cite[Thm.\,6.6]{Jon95}.  In particular, if $\alg{A}_p\notin\cl{X}$, then $\alg{A}_p\notin\Vop(\cl{X})$, and if $\alg{A}_{p-1}^+\notin\cl{X}$, then $\alg{A}_{p-1}^+\notin\Vop(\cl{X})$.  We have therefore established the following:
\begin{corollary}
There are $2^{\aleph_0}$ distinct varieties of semilinear Dunn (and likewise De Morgan) monoids without the weak ES property.
\end{corollary}

\section{Semilinear De Morgan Monoids}
\label{sec:DMM}

Except for the consequences of Theorem~\ref{thm:Aps}, we were concerned in Section~\ref{sec:no-involution} with involution-less algebras.
We now focus on algebras \emph{with} involution, and on varieties of De Morgan monoids.
Negatively generated totally ordered De Morgan monoids need not be idempotent (unlike their non-involutive counterparts).  We shall prove a representation theorem for these algebras, which will allow us to show that the negatively generated semilinear De Morgan monoids form a locally finite variety, all of whose subvarieties have the ES property.  The following lemma is well known.

\begin{lemma}\label{lem:bounds}
\textup{(\cite[Lem.~2.3]{MRW})}
If a (possibly involutive) RL\/ $\alg{A}$ has a least element\/ $\bot$\textup{,} then\/
$\top\seteq\bot\to\bot$ is its greatest element
and, for all\/
$a\in A$\textup{,}
\[
a\bdot\bot=\bot=\top\to\bot
\text{ \ and\/ \ }
\bot\to a=\top=
a\to\top=\top^2.
\]
In particular, $\{\bot,\top\}$ is a subalgebra of the\/ $\bdot,\to,\wedge,\vee\,(,\neg)$
reduct of\/ $\alg{A}$\textup{.}
\end{lemma}

If we say that $\bot,\top$ are \emph{extrema}
of an [I]RL $\alg{A}$, we mean that $\bot\leqslant a\leqslant\top$
for all $a\in A$.  An [I]RL with extrema is said to be \emph{bounded}.
In that case, its extrema need not be
\emph{distinguished} elements, and they are not always retained in subalgebras (consider the Sugihara monoids $\alg{S}_n$, for instance).
The next lemma is a straightforward consequence
of (\ref{eq:residuation}).
\begin{lemma}\label{lem:RC-conditions}
The following conditions on a bounded IRL $\alg{A}$\textup{,} with extrema\/ $\bot,\top$\textup{,}
are equivalent.
\begin{enumerate}
\item
$\top\bdot a=\top$ whenever\/ $\bot\neq a\in A$\textup{.}
\item
$a\to\bot=\bot$ whenever\/ $\bot\neq a\in A$\textup{.}
\item
$\top\to b=\bot$ whenever\/ $\top\neq b\in A$\textup{.}
\end{enumerate}
\end{lemma}
Following Meyer \cite{Mey86}, we say that an
IRL is
\emph{rigorously compact}
if it is bounded and satisfies
the equivalent conditions of Lemma~\ref{lem:RC-conditions}.
The next theorem is proved in \cite[Thm.~5.3]{MRW}, but has an antecedent in \cite[Thm.~3]{Mey86}.
\begin{theorem}
\label{thm:bounded-FSI-DMM=>RC}
Every bounded FSI De Morgan monoid is rigorously compact.
\end{theorem}

We depict below the two-element Boolean algebra $\mathbf{2}$, and two 
four-element
De Morgan monoids, $\alg{C}_4$ and $\alg{D}_4$.
In each case, the labeled Hasse diagram determines the structure.

{\tiny

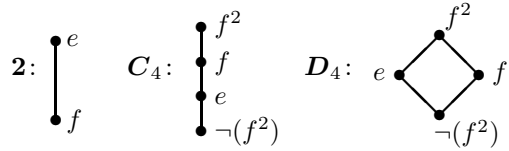
\thicklines
\begin{center}
\begin{picture}(80,60)(2,51)

\put(-25,63){\line(0,1){30}}
\put(-25,63){\circle*{4}}
\put(-25,93){\circle*{4}}

\put(-21,91){\small $e$}
\put(-21,60){\small $f$}

\put(-42,80){\small $\mathbf{{2}\colon}$}

%%%%%%

\put(30,59){\circle*{4}}
\put(30,59){\line(0,1){39}}
\put(30,72){\circle*{4}}
\put(30,85){\circle*{4}}
\put(30,98){\circle*{4}}

\put(35,96){\small ${f^2}$}
\put(35,82){\small $f$}
\put(35,69){\small ${e}$}
\put(35,56){\small $\neg(f^2)$}

\put(2,80){\small ${\alg{C}_4}\colon$}

%%%%%%

\put(120,65){\circle*{4}}
\put(135,80){\line(-1,-1){15}}
\put(135,80){\circle*{4}}
\put(105,80){\line(1,-1){15}}
\put(105,80){\circle*{4}}
\put(105,80){\line(1,1){15}}
\put(120,95){\circle*{4}}
\put(135,80){\line(-1,1){15}}

\put(122,99){\small ${f^2}$}
\put(95,78){\small ${e}$}
\put(140,78){\small $f$}
\put(118,55){\small $\neg(f^2)$}

\put(69,80){\small ${\alg{D}_4}\colon$}

\end{picture}\nopagebreak
\end{center}

}

\noindent
Note that a De Morgan monoid is $0$-generated iff it has no proper subalgebra.  The following result is implicit in Slaney \cite{Sla85,Sla89} and explicit in \cite[Thm.~5.20]{MRW}.

\begin{theorem}
\label{thm:0-gen-simples}
A De Morgan monoid is simple and $0$-generated iff it is isomorphic to $\alg{2}$ or to
$\alg{C}_4$ or to $\alg{D}_4$\textup.
\end{theorem}

Of these algebras, $\alg{C}_4$ garners special attention, because of the following. 
\begin{theorem}[{Slaney \cite[Thm.~1]{Sla89}}]\label{thm:Slaney}
Let $h\colon\alg{A} \rightarrow \alg{B}$ be a homomorphism, where $\alg{A}$ is an FSI De Morgan monoid, and $\alg{B}$ is nontrivial and\/ $0$-generated. Then $h$ is an isomorphism or $\alg{B} \cong \alg{C}_4$\textup.
\end{theorem}
A De Morgan monoid $\alg{A}$ is said to be \emph{crystalline}
if there is a homomorphism $h\colon\alg{A}\rightarrow\alg{C}_4$
(in which case $h$ is surjective).  These algebras do not form a variety, as their homomorphic images need not be crystalline, but
there is a largest variety $\cl{U}$ of crystalline (or trivial) De Morgan monoids; it is axiomatized in \cite[Sec.~4]{MRW2}.  Thus, $\Vop(\alg{C}_4)$ is the smallest nontrivial subvariety of $\cl{U}$.

\begin{theorem}[{\cite[Lem.~4.8]{MRW2}}]
\label{lem:W+RC=>U}
Let $\alg{A}$ be a rigorously compact and crystalline De Morgan monoid\textup.
Then $\alg{A} \in \cl{U}$\textup{.}
\end{theorem}

We say that a De Morgan monoid
is \emph{anti-idempotent} if it
satisfies $x\leqslant f^2$ (or equivalently,
$\neg(f^2)\leqslant x$).
This terminology is justified, because
a variety of square-increasing IRLs has no nontrivial idempotent member iff it satisfies $x\leqslant f^2$ \cite[Cor.~3.6]{MRW}.

We explained the structure of idempotent De Morgan monoids (i.e., Sugihara monoids) in Section~\ref{sec:preliminaries}, and recalled in Theorem~\ref{thm:SM-Es} that all varieties of Sugihara monoids have surjective epimorphisms. It is now convenient (in view of the upcoming Theorem~\ref{thm:FSI-DMMs}) to describe the anti-idempotent negatively generated totally ordered De Morgan monoids.

\begin{theorem}
\label{thm:DMM-neg-cone}
Let $\alg{A}$ be an anti-idempotent negatively generated FSI De Morgan monoid. 
Then $\alg{A} \cong \alg{D}_4$ or $\alg{A} \in \cl{U}$\textup.
\end{theorem}

\begin{proof}
We may suppose that $\alg{A}$ is nontrivial.
Being anti-idempotent, $\alg{A}$ has no trivial subalgebra, by Theorem~\ref{thm:idempotence}.
In a variety whose nontrivial members lack trivial subalgebras, every nontrivial member has a simple homomorphic image
\cite[Cor.~5.4]{MRW5}, so $\alg{A}$ has a simple homomorphic image $\alg{B}$.
Now $\alg{B} = \SgA^{\alg{B}} B^{-}$, by Lemma~\ref{lem:neg-cone-gen-in-images}, because
$\alg{A} = \SgA^{\alg{A}} A^{-}$.
Since $\alg{B}$ is simple and anti-idempotent, Lemma~\ref{lem:RL-properties}(\ref{lem:RL-properties:simple}) shows that
$B^{-}$ is the chain $\neg(f^2) < e$, 
so $\alg{B}$ is $0$-generated. Therefore, $\alg{B}$ is isomorphic to $\alg{C}_4$ or $\alg{D}_4$, by Theorem~\ref{thm:0-gen-simples}, as $\alg{2}$ is not anti-idempotent. If $\alg{B} \cong \alg{D}_4$, then $\alg{A} \cong \alg{D}_4$, by Theorem~\ref{thm:Slaney}.
Otherwise $\alg{B} \cong \alg{C}_4$, in which case 
$\alg{A}$ is crystalline (as well as FSI and bounded), so $\alg{A} \in \cl{U}$, by Theorems~\ref{thm:bounded-FSI-DMM=>RC} and \ref{lem:W+RC=>U}.
\end{proof}

The algebras in $\cl{U}$
are subdirect products of `skew reflections' of Dunn monoids
\cite[Cor.~5.6]{MRW2}.  The skew reflection
construction is a means of embedding
a Dunn monoid
into one that has an involution (i.e., into a De Morgan monoid).
In the semilinear context, to which we now confine ourselves, this
construction reduces to an older and simpler one, called
`reflection', which is recalled below.  It
is essentially due to Meyer \cite{Mey73}.

Given a Dunn monoid $\alg{A}$, let $A'=\{a':a\in A\}$ be a disjoint copy of $A$, and let $\0,\mathsf{\mathit{1}}$ be distinct non-elements of $A \cup A'$. The \emph{reflection} $\refl(\alg{A})$ of $\alg{A}$ is the De Morgan monoid with universe $\refl(A) = A \cup A' \cup \{\0,\mathsf{\mathit{1}}\}$ such that $\alg{A}$ is a subalgebra of the RL-reduct of $\refl(\alg{A})$ and, for all $a,b \in A$ and $x,y \in \refl(A)$,
\begin{align*}
& x \bdot \0 = \0 < a < b' < \mathsf{\mathit{1}} = a' \bdot b', \text{ and if } x \neq \0, \text{ then } x \bdot \mathsf{\mathit{1}} = \mathsf{\mathit{1}};\\
& a \bdot b' = (a \to b)';\\
& \neg a = a' \text{ and } \neg(a') = a \text{ and } \neg \0 = \mathsf{\mathit{1}} \text{ and } \neg \mathsf{\mathit{1}} = \0.
\end{align*}
Since $f = e'$, we have $\mathsf{\mathit{1}} = f^2$ and $\0 = \neg(f^2)$, 
so reflections are anti-idempotent.  Note that $\alg{C}_4\cong\refl(\alg{A})$ for any trivial Dunn monoid $\alg{A}$.

The \emph{reflection} of a variety $\cl{K}$ of Dunn monoids is the 
variety
\[ \Rop(\cl{K}) \seteq \Vop(\{\refl(\alg{A}) : \alg{A} \in \cl{K} \}). \]
We shall use the following facts concerning reflections, whose proofs can be found in \cite[Lem.~6.5]{MRW2} and \cite[Cor.~9.2, Thm.~9.3]{MRW4}:

\begin{theorem}
\label{thm:refl-props}
Let\/ $\cl{K}$ be a variety of Dunn monoids\textup. 
\begin{enumerate}
\item \label{thm:refl-props:S} If\/ $\alg{C}$ is a subalgebra of a Dunn monoid $\alg{D}$, then
\[
C \cup \{c' : c \in C\} \cup \{\mathsf{\mathit{1}},\0\}
\]
is the universe of a subalgebra of\/ $\refl(\alg{D})$ that is isomorphic to\/ $\refl(\alg{C})$, and every subalgebra $\alg{A}$ of\/ $\refl(\alg{D})$ arises in this way from a subalgebra $\alg{C}$ of $\alg{D}$, where $C = A \cap D = \{a \in A : a \neq \0 \text{ and } a^2 \neq \mathsf{\mathit{1}} \}$\textup. 

\item\label{lem:refl-props:H}
If\/ $\theta$ is a congruence of a Dunn monoid\/ $\alg{B}$\textup{,} then
\[
\quad\quad \refl(\theta)\seteq\theta\cup\{\langle a',b'\rangle : \langle a,b\rangle\in \theta\}\cup\{\langle \0,\0\rangle,\,\langle \mathsf{\mathit{1}},\mathsf{\mathit{1}}\rangle\}
\]
is a congruence of\/ $\refl(\alg{B})$\textup{,} and\/ $\refl(\alg{B})/\refl(\theta)\cong\,\refl(\alg{B}/\theta)$\textup{.}  Also, every proper congruence of\/
$\refl(\alg{B})$ has the form\/ $\refl(\theta)$ for some\/ $\theta\in \Con\,\alg{B}$\textup{.}

%\smallskip

\item\label{lem:refl-props:Pu}
If\/ $\{\alg{B}_i : i\in I\}$ is a family of Dunn monoids and\/ $\mathcal{U}$ is an ultrafilter over\/ $I$\textup{,} then\/
$\prod_{i\in I}\refl(\alg{B}_i)/\mathcal{U}\,\cong\refl\!\left(\prod_{i\in I}\alg{B}_i/\mathcal{U}\right)$\textup{.}

%\smallskip

\item \label{thm:refl-props:R-FSI} $\alg{A} \in \Rop(\cl{K})_\textup{FSI}$ iff\/ $\alg{A}$ is trivial or $\alg{A} \cong \refl(\alg{D})$ for some $\alg{D} \in \cl{K}_\textup{FSI}$\textup. 
    
\item \label{thm:refl-props:R-SI} $\alg{A} \in \Rop(\cl{K})_\textup{SI}$ iff\/
$\alg{A} \cong \refl(\alg{D})$, where $\alg{D}$ is trivial or belongs to $\cl{K}_\textup{SI}$\textup. 

\item \label{thm:refl-props:ES-refl} $\cl{K}$ has the ES property iff\/ $\Rop(\cl{K})$ has the ES property\textup. 

\item \label{thm:refl-props:refl-loc-fin} $\cl{K}$ is locally finite iff\/ $\Rop(\cl{K})$ is locally finite. More specifically, if $p \colon \omega \rightarrow \omega$ is a function such that, for each $n \in \omega$, every $n$-generated member of\/ $\cl{K}_\textup{FSI}$ has at most $p(n)$ elements, then every $n$-generated member of\/ $\Rop(\cl{K})_\textup{FSI}$ has at most\/ $2+2p(n)$ elements\textup. 
\end{enumerate}
\end{theorem}
In
(\ref{lem:refl-props:H}), if $\theta=B\times B$, then $\refl(\alg{B}/\theta)\cong\alg{C}_4$, so $\refl(\alg{B})\in\cl{U}$.

Recall that $\cl{S}$ is the class of all simple totally ordered idempotent RLs.
It follows
from Theorems~\ref{thm:SLid-simple-ES} and \ref{thm:refl-props}(\ref{thm:refl-props:ES-refl}) that $\Rop(\Vop(\cl{S}))$ has the ES property.
It also has finite depth and its members are not all negatively generated.

\begin{lemma}
\label{lem:SLaiDMM-neg-cone}
Every nontrivial totally ordered negatively generated anti-idempotent De Morgan monoid $\alg{A}$ is a reflection of a totally ordered Dunn monoid. 
\end{lemma}

\begin{proof}
As $\alg{A}$ is negatively generated, FSI and anti-idempotent, Theorem~\ref{thm:DMM-neg-cone} shows that $\alg{A} \in \cl{U}\cup\Iop(\alg{D}_4)$.
But $\alg{D}_4$ is not totally ordered, so $\alg{A} \in \cl{U}$.
Since $\alg{A}$ is bounded and FSI, it
is rigorously compact, by Theorem~\ref{thm:bounded-FSI-DMM=>RC}.
Also, $\alg{A}$ is crystalline, like every nontrivial member of $\cl{U}$. These two properties (being rigorously compact and crystalline) are enough to guarantee that $\alg{A}$ is a `skew reflection' of a Dunn monoid
$\alg{B}$, by \cite[Thm.~5.4]{MRW2}.
In the present context, since $\alg{A}$ is totally ordered, this amounts to saying that $\alg{A}$ is a
reflection of $\alg{B}$, which is also totally ordered.
\end{proof}

The underlying Dunn monoid in the statement of Lemma~\ref{lem:SLaiDMM-neg-cone}
is itself negatively generated, because of the next lemma.
We shall see in the proof of Theorem~\ref{thm:SLaiDMM-neg-cone} that the converse of Lemma~\ref{lem:SLaiDMM-neg-cone} holds for such (negatively generated) Dunn monoids.

\begin{lemma}
\label{lem:refl-subalg-gen}
Let $\alg{A} = \refl(\alg{D})$ for some Dunn monoid $\alg{D}$\textup.
If\/ $\alg{A} = \SgA^{\alg{A}} X$ for some $X \subseteq D$, then $\alg{D}=\SgA^{\alg{D}} X$\textup.
\end{lemma}

\begin{proof}
Let $\alg{B}$ be the subalgebra of $\alg{D}$ generated by $X$.
We argue that $\alg{B} = \alg{D}$.
By Lemma~\ref{thm:refl-props}(\ref{thm:refl-props:S}), $\refl(\alg{B})$ can be identified with a subalgebra of $\alg{A}$.
But then $\refl(\alg{B}) = \alg{A} = \refl(\alg{D})$, since $\alg{A} = \SgA^{\alg{A}} X$ and $X \subseteq B \subseteq \refl(B)$. It follows that $B = D$, because $\alg{A}$ is a reflection of a Dunn monoid whose universe must be
$\{a \in A : a \neq \0 \text{ and } a^2 \neq \mathsf{\mathit{1}}\}$ (again by Lemma~\ref{thm:refl-props}(\ref{thm:refl-props:S})).
\end{proof}

We can now prove a representation theorem for semilinear negatively generated anti-idempotent De Morgan monoids, which also reveals the unobvious fact that these algebras form a variety.
We define the following unary terms:
\begin{align*}
d'(x) &\seteq (f^2 \to (x \bdot f)) \wedge (f^2 \bdot \neg x);\\
\sigma(x) &\seteq (x \wedge e) \bdot (x^{*} \wedge e)^{*};\\
d(x) &\seteq d'(\neg x) \,\text{ and }\, \sigma'(x) \seteq \neg \sigma(\neg x).
\end{align*}
Recall that $\sigma(x) = x$ is equation (\ref{eq:sigma}), 
which is satisfied by every generalized Sugihara monoid. Consider the equation
\begin{equation}
\label{eq:negcone}
x = \big(d(\sigma(x)) \wedge \sigma(x)\big) \vee \big(d'(\sigma'(x)) \wedge \sigma'(x)\big) \vee \left(\big(f^2 \vee \neg(f^2)\big) \to \sigma'(x)\right)\!.
\end{equation}
Of course, $f^2 \vee \neg(f^2)$ amounts to $f^2$ in anti-idempotent De Morgan monoids.  We have not exploited this simplification in (\ref{eq:negcone}), because
the next result will be generalized in Theorem~\ref{thm:SL-DMM-neg-cone} to accommodate De Morgan monoids that need \emph{not} be anti-idempotent. 

\begin{theorem}
\label{thm:SLaiDMM-neg-cone}
Let $\alg{A}$ be a nontrivial anti-idempotent semilinear De Morgan monoid. Then the following are equivalent:
\begin{enumerate}
\item \label{thm:SLaiDMM-neg-cone:neg-cone} $\alg{A}$ is negatively generated;
\item \label{thm:SLaiDMM-neg-cone:FSI} $\alg{A}$ is a subdirect product of reflections of totally ordered generalized Sugihara monoids;
\item \label{thm:SLaiDMM-neg-cone:equation} $\alg{A}$ satisfies equation~\textup{(\ref{eq:negcone}).}
\end{enumerate}
\end{theorem}

\begin{proof}
(\ref{thm:SLaiDMM-neg-cone:neg-cone}) $\Rightarrow$ (\ref{thm:SLaiDMM-neg-cone:FSI}):
As in the proof Theorem~\ref{thm:SLDunn-neg-cone}, it suffices, by the
Subdirect Decomposition Theorem, to show that every nontrivial \emph{totally ordered} anti-idempotent De Morgan monoid $\alg{B}$ that is negatively generated is a reflection of a totally ordered generalized Sugihara monoid. By Lemma~\ref{lem:SLaiDMM-neg-cone}, $\alg{B} \cong \refl(\alg{D})$ for some totally ordered Dunn monoid $\alg{D}$. Note that $\refl(\alg{D})$ is generated by $D^{-}$, because $\refl(D)^{-} = D^{-} \cup \{\0\}$ and $\0 = \neg(f^2) \in \Sg^{\alg{B}}\{e\}$. But then $\alg{D} = \SgA^{\alg{D}} D^{-}$, by Lemma~\ref{lem:refl-subalg-gen}. It follows, by Theorem~\ref{thm:SLDunn-neg-cone}, that $\alg{D} \in \cl{GSM}$. 

(\ref{thm:SLaiDMM-neg-cone:FSI}) $\Rightarrow$ (\ref{thm:SLaiDMM-neg-cone:equation}): We claim that every reflection of a totally ordered generalized Sugihara monoid satisfies (\ref{eq:negcone}), and so $\alg{A}$ does as well. Let $\alg{B} = \refl(\alg{D})$ for some totally ordered $\alg{D} \in \cl{GSM}$.
For any $a \in B$, it follows from the definition of reflection that
\[
d(a) = \begin{cases}
\mathsf{\mathit{1}} &\text{if } a \in D; \\
\0 &\text{otherwise},
\end{cases}
\quad
d'(a) = \begin{cases}
\mathsf{\mathit{1}} &\text{if } a \in D'; \\
\0 &\text{otherwise},
\end{cases}
\]
\[
(f^2 \vee \neg(f^2)) \to a =
f^2 \to a =
\begin{cases}
\mathsf{\mathit{1}} &\text{if } a = \mathsf{\mathit{1}}; \\
\0 &\text{otherwise},
\end{cases}
\]
\[
\sigma(a) = \begin{cases}
\mathsf{\mathit{1}} &\text{if } a \in D'; \\
a &\text{otherwise},
\end{cases}
\;\text{ and }\;
\sigma'(a) = \begin{cases}
\0 &\text{if } a \in D; \\
a &\text{otherwise}.
\end{cases}
\]
It is then easy to verify that $\alg{B}$ satisfies (\ref{eq:negcone}), by checking the cases where $a \in D$, $a \in D'$, $a = \mathsf{\mathit{1}}$ and $a = \0$.

(\ref{thm:SLaiDMM-neg-cone:equation}) $\Rightarrow$ (\ref{thm:SLaiDMM-neg-cone:neg-cone}): This follows directly from the shape of equation (\ref{eq:negcone}),
because $\sigma$ is built up from the terms $x \wedge e$ and $x^{*} \wedge e$,
and
$\sigma'$ is built up from
$\neg x \wedge e$ and $(\neg x)^{*} \wedge e$. For any assignment to $x$ of an element of $\alg{A}$, these terms clearly evaluate into $A^{-}$.
\end{proof}

\begin{corollary}
\label{cor:SlaiDMM-loc-fin}
Let\/ $\cl{K}$ be the class of negatively generated semilinear anti-idempotent De Morgan monoids. Then
\begin{enumerate}
\item\label{cor:SlaiDMM-loc-fin:variety} $\cl{K}$ is a variety that is axiomatized relative to semilinear De Morgan monoids by $x \leqslant f^2$ and\/ \textup{(\ref{eq:negcone});}
\item\label{cor:SlaiDMM-loc-fin:Rop} $\cl{K} = \Rop(\cl{GSM})$\textup;
\item\label{cor:SlaiDMM-loc-fin:bound} if\/ $\alg{A} \in \cl{K}$ is totally ordered and $n$-generated, then\/ $|A| \leq 6n + 4$\textup;
\item\label{cor:SlaiDMM-loc-fin:loc-fin} $\cl{K}$ is locally finite.
\end{enumerate}
\end{corollary}
\begin{proof}
(\ref{cor:SlaiDMM-loc-fin:variety}) follows immediately from Theorem~\ref{thm:SLaiDMM-neg-cone}.

For (\ref{cor:SlaiDMM-loc-fin:Rop}), it follows straightforwardly from Theorem~\ref{thm:SLaiDMM-neg-cone} that $\cl{K} \subseteq \Rop(\cl{GSM})$.
To establish the converse, it is enough to show that $\Rop(\cl{GSM})_\textup{SI} \subseteq \cl{K}$, because $\cl{K}$ is closed under $\Iop\Psop$ (by (\ref{cor:SlaiDMM-loc-fin:variety})). By Theorem~\ref{thm:refl-props}(\ref{thm:refl-props:R-SI}), this reduces to showing that $\Rop(\cl{GSM}_\textup{SI}) \subseteq \cl{K}$, which follows from Theorem~\ref{thm:SLaiDMM-neg-cone}.

Recall from 
the remarks after Lemma~\ref{fct:FG-function}
that
if $\alg{B} \in \cl{GSM}$ is totally ordered and $n$-generated then $|B| \leq 3n+1$.
Let $\alg{A}$ be a totally ordered $n$-generated member of $\cl{K}$.
By Theorem~\ref{thm:refl-props}(\ref{thm:refl-props:refl-loc-fin}), 
$|A|\leq 2+2(3n+1) = 6n + 4$, proving (\ref{cor:SlaiDMM-loc-fin:bound}).

(\ref{cor:SlaiDMM-loc-fin:loc-fin}) follows from
(\ref{cor:SlaiDMM-loc-fin:bound}) (and Lemma~\ref{fct:FG-function}).
\end{proof}

\begin{corollary}
\label{cor:SlaiDMM-subv-Rop}
Let $\cl{K}$ be any nontrivial variety of negatively generated semilinear anti-idempotent De Morgan monoids. Then $\cl{K} = \Rop(\cl{L})$ for some variety\/ $\cl{L}$ of generalized Sugihara monoids.
\end{corollary}
\begin{proof}
Let $\cl{D} =
\{\alg{D} \in \cl{GSM} :  \refl(\alg{D}) \in \cl{K}_\textup{SI}\}$
and
$\cl{L} = \Vop(\cl{D})$. By the
Subdirect Decomposition Theorem, it suffices to show that $\cl{K}_\textup{SI} = \Rop(\cl{L})_\textup{SI}$.

Let $\alg{A} \in \cl{K}_\textup{SI}$. By (\ref{thm:SLaiDMM-neg-cone:neg-cone}) $\Rightarrow$ (\ref{thm:SLaiDMM-neg-cone:FSI}) of Theorem~\ref{thm:SLaiDMM-neg-cone}, $\alg{A} \cong\refl(\alg{D})$ for some $\alg{D} \in \cl{GSM}$. But then $\alg{D} \in \cl{D}$, so $\alg{A} \cong \refl(\alg{D}) \in \Rop(\cl{L})$.

Conversely, let $\alg{A} \in \Rop(\cl{L})_\textup{SI}$. By Theorem~\ref{thm:refl-props}(\ref{thm:refl-props:R-SI}), $\alg{A} \cong \refl(\alg{D})$ for some $\alg{D} \in \cl{L}$ that is either trivial or SI. In the first case $\alg{A} \cong \alg{C}_4$, and $\alg{C}_4 \in \cl{K}$, because $\cl{K}$ is a nontrivial subvariety of $\cl{U}$ 
(Lemma~\ref{lem:SLaiDMM-neg-cone}). In the second case, $\alg{D} \in \Vop(\cl{D})_\textup{SI} \subseteq \Hop\Sop\Puop(\cl{D})$, by J\'{o}nsson's Theorem. So, by Lemma~\ref{thm:refl-props}(\ref{thm:refl-props:S})--(\ref{lem:refl-props:Pu}), 
\[\alg{A} \cong \refl(\alg{D}) \in \Hop\Sop\Puop(\{\refl(\alg{B}):\alg{B} \in \cl{D}\}) \subseteq \cl{K}.
\qedhere
\]
\end{proof}

\begin{theorem}
\label{thm:SLaiDMM-ES}
Let\/ $\cl{K}$ be any variety of negatively generated semilinear anti-idempotent De Morgan monoids.
Then\/ $\cl{K}$ has surjective epimorphisms.
\end{theorem}

\begin{proof}
We may suppose without loss of generality that $\cl{K}$ is nontrivial, so ${\cl{K} = \Rop(\cl{L})}$ for some variety $\cl{L}$ of generalized Sugihara monoids, by Corollary~\ref{cor:SlaiDMM-subv-Rop}.
By Theorem~\ref{thm:GSM-ES}, $\cl{L}$ has surjective epimorphisms, so by
Theorem~\ref{thm:refl-props}(\ref{thm:refl-props:ES-refl}), $\Rop(\cl{L}) = \cl{K}$ has as well.
\end{proof}

We aim now to generalize the above results by dropping anti-idempotency, so our focus will be on
negatively generated semilinear De Morgan monoids in general.
We first recall some structural facts about De Morgan monoids.  As usual, in a poset, we denote by $(a]$ the set of all lower bounds of an element $a$ (including $a$ itself), and by $[a)$ the set of all upper bounds.
\begin{theorem}[{\cite[Thm.~5.15--18]{MRW}}]
\label{thm:non-idempotentDMMs}
Let $\alg{A}$ be a non-idempotent FSI De Morgan monoid.
\begin{enumerate}
\item \label{thm:interval-subalgebra} 
If $f^2\leqslant a\in A$, then $\neg a<a$ and the interval
$[\neg a,a]$ is a subuniverse of $\alg{A}$.

\item \label{thm:lollipop} 
$\alg{A}$ is the union of the interval subuniverse $[\neg(f^2),f^2]$ and two chains of idempotents,
$(\neg(f^2)]$ and\/ $[f^2)$.

\item \label{thm:fusion-behavior} 
If $f^2 \leqslant a < b$, then $a \to a = a$, $a \to b = b$ and $b \to a = \neg b$.
\end{enumerate}
\end{theorem}

It follows from (\ref{thm:interval-subalgebra}) that $\neg(f^2)\leqslant e$, so $[\neg(f^2))$ is a deductive filter of $\alg{A}$.

\begin{theorem}\label{thm:Sugihara-factor-odd}
Let\/ $\alg{A}$ be a non-idempotent FSI De Morgan monoid. Then\/ $\alg{A}/[\neg(f^2))$ is a totally ordered odd Sugihara monoid.
Furthermore, $e/[\neg(f^2))$ is the interval $[\neg(f^2), f^2]$, and $a/[\neg(f^2)) = \{a\}$ for any $a \in A \setminus [\neg(f^2), f^2]$.
\end{theorem}
\begin{proof}
Let $G \seteq [\neg(f^2))$
and $a \in [\neg(f^2), f^2]$.
By Theorem~\ref{thm:non-idempotentDMMs}(\ref{thm:interval-subalgebra}), $[\neg(f^2), f^2]$ is a subuniverse of $\alg{A}$, so $e \to a,\; a \to e \in [\neg(f^2), f^2] \subseteq G$, whence $a/G = e/G$.
Therefore, $[\neg(f^2), f^2] \subseteq e/G$.
In particular, since $f \in [\neg(f^2), f^2]$, we have $e/G = f/G$, so $\alg{A}/G$ is an odd Sugihara monoid, by Theorem~\ref{thm:odd-DMMs-are-SM}.
By Theorem~\ref{thm:non-idempotentDMMs}(\ref{thm:lollipop}), $A \setminus [\neg(f^2), f^2]$ is totally ordered, so $\alg{A}/G$ is as well.

Let $a \in e/G$. Then $\neg(f^2) \leqslant a$ and $\neg(f^2) \leqslant a \to e$. By the law of residuation $a \bdot \neg(f^2) \leqslant e$, so by (\ref{eq:neg-fusion-law}), $\neg(f^2) \bdot f \leqslant \neg a$.
Since $[\neg(f^2), f^2]$ is a subuniverse of $\alg{A}$ with least element $\neg(f^2)$, we have $\neg(f^2) = \neg(f^2) \bdot f \leqslant \neg a$, by Lemma~\ref{lem:bounds}. So, $a \leqslant f^2$. Therefore $e/G = [\neg(f^2), f^2]$.

Lastly, let $a \in A \setminus [\neg(f^2), f^2]$, and suppose that $a/G = b/G$ for some $b \in A$. Notice that $b \notin [\neg(f^2), f^2]$, since $a \notin e/G = [\neg(f^2), f^2]$.

By involutional symmetry, 
we may assume  that $f^2 < a$ (rather than $a < \neg(f^2)$), because otherwise $f^2 < \neg a$, and from $x/G = \{x\}$ and the double negation law, it follows easily that $(\neg x)/G = \{\neg x\}$.

If $b < \neg(f^2)$, then $b < e < a$, but $a/G$ includes $a$ and $b$, and is an interval of $\alg{A}$, so it includes $e$, whence $a/G = e/G$, a contradiction.
Therefore, $f^2 < b$.
By Theorem~\ref{thm:non-idempotentDMMs}(\ref{thm:fusion-behavior}),
\[a \to b  \in \{a,b,\neg a, \neg b\} \subseteq A \setminus [\neg(f^2), f^2].\]
As $a/G = b/G$, we have $\neg(f^2) \leqslant a \to b,\, b \to a$,
so $e < f^2 < a \to b$. Similarly, $e < b \to a$, so $a = b$. Therefore, $a/G = \{a\}$.
\end{proof}

The following discussion elaborates and systematizes
Remark~5.19 of \cite{MRW},
by showing
how any non-idempotent FSI De Morgan monoid can be viewed as an
extension of its anti-idempotent subalgebra on $[\neg(f^2), f^2]$ by the (idempotent) totally ordered odd Sugihara monoid that results from factoring out $[\neg(f^2))$.  We call this a `rigorous extension', as it is a union of rigorously compact algebras containing $[\neg(f^2),f^2]$.

Let $\alg{S}$ be a totally ordered odd Sugihara monoid. 
For any non-constant basic operation $\varphi$ of $\alg{S}$ with arity $n > 0$, and for any $a_1, \dots, a_n \in S$,
\begin{equation}
\label{eq:OSM-totally-irreducible}
\text{if } \varphi(a_1, \dots, a_n) = e \text{ then } a_i = e \text{ for some } i \leq n.
\end{equation}
When $\varphi$ is $\neg$, (\ref{eq:OSM-totally-irreducible}) follows from the fact that $\alg{S}$ is odd, and when $\varphi$ is $\wedge$ or $\vee$, (\ref{eq:OSM-totally-irreducible}) holds because $\alg{S}$ is totally ordered. When $\varphi$ is $\bdot$, notice that the odd Sugihara monoid $\alg{Z}$ satisfies the quasi-equation $x \bdot y = e \implies x = e$, so since $\cl{OSM}$ is generated as a quasivariety by $\alg{Z}$, 
$\alg{S}$ satisfies the same quasi-equation, whence (\ref{eq:OSM-totally-irreducible}) holds.

Except for the treatment of involution, the construction in the next definition coincides 
with one in Galatos \cite[p.\,458]{Gal11}.
\begin{definition}\label{def:rigorous-extension}
The \emph{rigorous extension}
of a De Morgan monoid $\alg{A}$ by a totally ordered odd Sugihara monoid $\alg{S}$
is the algebra
\[\alg{S}[\alg{A}] \seteq \langle (S \setminus \{e^{\alg{S}}\}) \cup A; \bdot', \wedge', \vee', \neg', e^{\alg{A}} \rangle\]
with the following properties. Let $ \star \in \{\wedge,\vee,\bdot\}$.
The operations $\neg'$ and $\star'$ extend those of $\alg{S}$ and $\alg{A}$, i.e.,
for every $s,p \in S \setminus \{e^{\alg{S}}\}$ and $a,b \in A$, 
\[
\neg' s \seteq \neg^{\alg{S}} s,\;\;
\neg' a \seteq \neg^{\alg{A}} a,\;\;
s \star' p \seteq  s \star^{\alg{S}} p,\;\;\text{ and }\;\;
a \star' b \seteq  a \star^{\alg{A}} b
\]
(whence $\{\neg' s,\; s \star' p\} \subseteq S \setminus \{e^{\alg{S}}\}$, by (\ref{eq:OSM-totally-irreducible})), while
\[ a \star' s \seteq s \star' a \seteq \begin{cases}
a &\text{if } e^{\alg{S}} \star^{\alg{S}} s = e^{\alg{S}}; \\
e^{\alg{S}} \star^{\alg{S}} s &\text{otherwise.} 
\end{cases} \]
\end{definition}

\begin{theorem}
For any De Morgan monoid $\alg{A}$ and any totally ordered odd Sugihara monoid $\alg{S}$, the algebra $\alg{S}[\alg{A}]$ is a De Morgan monoid having $\alg{A}$ as a subalgebra.
\end{theorem}

\begin{proof}
It is easy to see that $\langle (S \setminus \{e^{\alg{S}}\}) \cup A; \wedge', \vee' \rangle$ is a lattice, that its lattice order $\leqslant$ extends $\leqslant^{\alg{A}}$ and $\leqslant^{\alg{S}}|_{S\setminus\{e^{\alg{S}}\}}$, and that for all $s \in S\setminus\{e^{\alg{S}}\}$ and $a \in A$ we have
\[(a \leqslant s \;\text{ iff }\; e^{\alg{S}} \leqslant^{\alg{S}} s) \;\text{ and }\; (s \leqslant a \;\text{ iff }\; s \leqslant^{\alg{S}} e^{\alg{S}}).\]
Since $\alg{S}$ is totally ordered and $\alg{A}$ distributive, the construction precludes diamond or pentagon sublattices, so
$\leqslant$ is a distributive lattice order.

It is straightforward to verify that $\bdot'$ is associative and has identity $e^\alg{A}$, and that (\ref{eq:neg-fusion-law}) is satisfied.
Here, it is helpful to note that there is no element $s \in S \setminus \{e^{\alg{S}}\}$ such that
$e^{\alg{S}} \bdot^{\alg{S}} s = e^{\alg{S}}$.
So, $s \bdot' a = a \bdot' s = s$ for every $s \in S \setminus \{e^{\alg{S}}\}$ and $a \in A$.
\end{proof}

\begin{theorem}
\label{thm:FSI-DMMs}
If\/ $\alg{A}$ is an FSI De Morgan monoid, then one of the following mutually exclusive conditions holds:
\begin{enumerate}
\item \label{thm:FSI-DMMs:SM} $\alg{A}$ is a Sugihara monoid, or
\item \label{thm:FSI-DMMs:AI} $\alg{A} \cong \alg{S}[\alg{A}']$, where $\alg{A}'$ is the nontrivial anti-idempotent subalgebra of $\alg{A}$ with universe $[\neg(f^2),f^2]$, and $\alg{S}$ is the totally ordered odd Sugihara monoid $\alg{A}/[\neg(f^2))$\textup.
\end{enumerate}
\end{theorem}

\begin{proof}
Let $\alg{A}$ be an FSI De Morgan monoid in which (\ref{thm:FSI-DMMs:SM}) fails. Then $\alg{A}$ is non-idempotent, with $f < f^2$. Let $G = [\neg(f^2))$ and $\alg{S} = \alg{A}/G$. Then $\alg{S}$ is a totally ordered odd Sugihara monoid, by Theorem~\ref{thm:Sugihara-factor-odd}.
Let $\alg{A}'$ be the
nontrivial anti-idempotent
subalgebra of $\alg{A}$ with universe $[\neg(f^2), f^2]$, which exists by Theorem~\ref{thm:non-idempotentDMMs}(\ref{thm:interval-subalgebra}).
We show that $\alg{A} \cong \alg{S}[\alg{A}']$, the isomorphism being 
\[
h \colon a \mapsto \begin{cases}
a &\text{ if } a \in A'; \\
a/G &\text{ otherwise.}
\end{cases}
\]
It follows from Theorem~\ref{thm:Sugihara-factor-odd} that $h$ is a bijection. It remains to show that $h$ is a homomorphism. It is clear that $h$ preserves $e$ and $\neg$. Let $\star \in \{\wedge,\vee,\bdot\}$. If $a,b \in A'$ then $h(a) \star^{\alg{S}[\alg{A}']} h(b) = a \star^{\alg{A}'} b = h(a \star^{\alg{A}} b)$, since $\alg{A}'$ is a subalgebra of $\alg{A}$ and of $\alg{S}[\alg{A}']$.
If $a,b \in A \setminus A'$, then $a \star^{\alg{A}} b \notin A'$, 
because otherwise $a/G \mathop{\star^{\alg{S}}} b/G = e/G$, whence $a/G = e/G$ or $b/G = e/G$, by (\ref{eq:OSM-totally-irreducible}), contradicting the fact that $a/G=\{a\}$ and $b/G=\{b\}$ (Theorem~\ref{thm:Sugihara-factor-odd}).
So,
\[
h(a) \star^{\alg{S}[\alg{A}']} h(b)
= a/G \star^{\alg{S}[\alg{A}']} b/G
= a/G \star^{\alg{S}} b/G
= (a \star^{\alg{A}} b)/G
= h(a \star^{\alg{A}} b).
\]
Now let $a \in A'$ and $b \in A \setminus A'$. If $e/G \wedge^{\alg{S}} b/G = e/G$ then $f^2 < b$, by Theorems~\ref{thm:non-idempotentDMMs}(\ref{thm:lollipop}) and~\ref{thm:Sugihara-factor-odd},
so
$
h(a) \wedge^{\alg{S}[\alg{A}']} h(b)
= a
= h(a \wedge^{\alg{A}} b)
$. 
If $e/G \wedge^{\alg{S}} b/G \neq e/G$ then $e/G \wedge^{\alg{S}} b/G = b/G$, as $\alg{S}$ is totally ordered. Then $b < \neg(f^2)$, so
$
h(a) \wedge^{\alg{S}[\alg{A}']} h(b)
= b/G
= h(a \wedge^{\alg{A}} b)
$. 
Similarly, $h(a) \vee^{\alg{S}[\alg{A}']} h(b) = h(a \vee^{\alg{A}} b)$.

It remains to show that $h(a) \bdot^{\alg{S}[\alg{A}']} h(b) = h(a \bdot^{\alg{A}} b)$. Note that
\[
h(a) \bdot^{\alg{S}[\alg{A}']} h(b) = e/G \bdot^{\alg{S}} b/G = b/G,
\]
so we must show that $a \bdot^{\alg{A}} b = b$.
This follows from the fact that $a$ and $b$ belong to the rigorously compact interval subalgebra of $\alg{A}$ with idempotent extrema $b$ and $\neg b$; see Theorems~\ref{thm:bounded-FSI-DMM=>RC} and \ref{thm:non-idempotentDMMs}(\ref{thm:interval-subalgebra}).
\end{proof}

Theorem~\ref{thm:FSI-DMMs} largely reduces the study of irreducible De Morgan monoids to the anti-idempotent case, about which we already have much information in the semilinear subcase.
The following properties of rigorous extensions are useful.

\begin{theorem}
\label{thm:expandingConstruction}
Let $\{\alg{A}, \alg{B}\} \cup \{\alg{A}_i : i \in I\}$ be a family of De Morgan monoids, and $\{\alg{S}\} \cup \{\alg{S}_i : i \in I\}$ a family of totally ordered odd Sugihara monoids, for some set $I$\textup.
\begin{enumerate}
\item \label{thm:expandingConstruction:H} If\/ $h \colon \alg{A} \rightarrow \alg{B}$ is a homomorphism, then the map
\[
h' \colon x \mapsto \begin{cases}
h(x) &\text{if\/ } x \in A; \\
x &\text{otherwise},
\end{cases}
\]
is a homomorphism from $\alg{S}[\alg{A}]$ to $\alg{S}[\alg{B}]$ which extends $h$\textup.
\item \label{thm:expandingConstruction:S} If\/ $\alg{P}$ is a subalgebra of $\alg{S}$ and $\alg{B}$ a subalgebra of $\alg{A}$, then $\alg{P}[\alg{B}]$ is a subalgebra of $\alg{S}[\alg{A}]$\textup.
\item \label{thm:expandingConstruction:Pu} $\prod_{i \in I} (\alg{S}_i[\alg{A}_i])/\mathcal{U} \cong \left(\prod_{i \in I} \alg{S}_i/\mathcal{U} \right) \left[\prod_{i \in I} \alg{A}_i/\mathcal{U}\right]$ for every ultrafilter $\mathcal{U}$ over $I$\textup.
\end{enumerate}
\end{theorem}

\begin{proof}
For (\ref{thm:expandingConstruction:H}),
we only show preservation of the binary basic operations with mixed arguments from $S\setminus\{e^{\alg{S}}\}$ and $A$, since the other cases are trivial.
Let $s \in S\setminus\{e^{\alg{S}}\}$ and $a \in A$.
If $s < a$ then $h'(s \wedge a) = h'(s) = s =h'(s) \wedge h'(a)$ and $h'(s \vee a) = h'(a) =h(a) = s \vee h(a) = h'(s) \vee h'(a)$.
When $a < s$, the argument is symmetrical. Also,
\begin{align*}
&h'(s \bdot a) = h'(s) = s = s \bdot h(a) = h'(s) \bdot h'(a).
\end{align*}

Item (\ref{thm:expandingConstruction:S}) follows from the fact that if $p \in P$ and $b \in B$, then for any $\star \in \{\wedge, \vee, \bdot\}$ we have $\{\neg p, \neg b, p \star b, b \star p\} \subseteq \{b, \neg b,p,\neg p\} \subseteq P[B]$.

In (\ref{thm:expandingConstruction:Pu}), we use the notation $\vec{x}=\langle x_i:i\in I\rangle$ for elements of
$\prod_{i \in I} S_i[A_i]$.
For $\vec{a}
\in \prod_{i \in I} S_i[A_i]$, let $I_{\vec{a}} \seteq \{i \in I : a_i \in A_i\}$. When $I_{\vec{a}} \in \mathcal{U}$, let $h(\vec{a}) = \vec{b}/\mathcal{U} \in \prod_{i \in I} A_i/\mathcal{U}$ where
\[b_i = a_i \text{ if } a_i \in A_i \text{, and } b_i = e^{\alg{A}_i} \text{ otherwise}.\]
When $I_{\vec{a}} \notin \mathcal{U}$, then its complement $I_{\vec{a}}^{c} = \{i \in I : a_i \in S_i\setminus\{e^{\alg{S}_i}\}\} \in \mathcal{U}$, since $\mathcal{U}$ is an ultrafilter. In this case, let $h(\vec{a}) = \vec{s}/\mathcal{U} \in (\prod_{i \in I} S_i/\mathcal{U}) \setminus \{e\}$ where
\[s_i = a_i \text{ if } a_i \in S_i \text{, and } s_i = e^{\alg{S}_i} \text{ otherwise}.\]
It can be verified that $h$ is a surjective homomorphism from $\prod_{i \in I} \alg{S}_i[\alg{A}_i]$ to $\left(\prod_{i \in I} \alg{S}_i/\mathcal{U} \right) \left[\prod_{i \in I} \alg{A}_i/\mathcal{U}\right]$,
whose
kernel is the congruence of $\prod_{i \in I} \alg{S}_i[\alg{A}_i]$ associated with $\mathcal{U}$. Then $\prod_{i \in I} (\alg{S}_i[\alg{A}_i])/\mathcal{U} \cong \left(\prod_{i \in I} \alg{S}_i/\mathcal{U} \right) \left[\prod_{i \in I} \alg{A}_i/\mathcal{U}\right]$, by the Homomorphism Theorem.
\end{proof}

\begin{corollary}
\label{cor:expandingHSPu}
Let $\alg{A}$ be a De Morgan monoid and $\alg{S}$ a totally ordered odd Sugihara monoid.
If\/ $\alg{C} \in \Hop \Sop \Puop (\alg{A})$, then $\alg{S}[\alg{C}] \in \Hop \Sop \Puop (\alg{S}[\alg{A}])$\textup.
\end{corollary}

\begin{proof}
Suppose $h \colon \alg{B} \rightarrow \alg{C}$ is a surjective homomorphism, with $\alg{B}$ a subalgebra of $\prod_{i \in I} \alg{A} /\mathcal{U}$ for some ultrafilter $\mathcal{U}$ over a set $I$.
By Theorem~\ref{thm:expandingConstruction}(\ref{thm:expandingConstruction:H}), $h$ can be extended to a surjective homomorphism $h'\colon\alg{S}[\alg{B}] \to \alg{S}[\alg{C}]$.
Recall that any algebra embeds into each of its ultrapowers. In particular, we may identify $\alg{S}$ with a subalgebra of $\prod_{i \in I} \alg{S} /\mathcal{U}$.
Then, by Theorem~\ref{thm:expandingConstruction}(\ref{thm:expandingConstruction:S}), $\alg{S}[\alg{B}]$ is a subalgebra of $\left(\prod_{i \in I} \alg{S} /\mathcal{U}\right) [\prod_{i \in I} \alg{A} /\mathcal{U}]$.
Lastly, by Theorem~\ref{thm:expandingConstruction}(\ref{thm:expandingConstruction:Pu}), $\left(\prod_{i \in I} \alg{S} /\mathcal{U}\right)[\prod_{i \in I} \alg{A} /\mathcal{U}] \cong \prod_{i \in I} (\alg{S}[\alg{A}]) /\mathcal{U}$.
So, $\alg{S}[\alg{C}] \in \Hop\Sop\Puop(\alg{S}[\alg{A}])$.
\end{proof}

We can now describe all semilinear De Morgan monoids that are negatively generated, using the characterization of FSI De Morgan monoids (in Theorem~\ref{thm:FSI-DMMs}) by means of rigorous extensions.

\begin{theorem}
\label{thm:SL-DMM-neg-cone}
Let $\alg{A}$ be a semilinear De Morgan monoid. Then the following are equivalent:
\begin{enumerate}
\item \label{thm:SL-DMM-neg-cone:neg-cone} $\alg{A}$ is negatively generated;
\item \label{thm:SL-DMM-neg-cone:FSI} $\alg{A}$ is a subdirect product of totally ordered Sugihara monoids and De Morgan monoids of the form $\alg{S}[\refl(\alg{D})]$\textup, where $\alg{S} \in \cl{OSM}_\textup{FSI}$ and $\alg{D} \in \cl{GSM}_\textup{FSI}$\textup;
\item \label{thm:SL-DMM-neg-cone:equation} $\alg{A}$ satisfies equation~\textup{(\ref{eq:negcone}).}
\end{enumerate}
\end{theorem}

\begin{proof}
(\ref{thm:SL-DMM-neg-cone:neg-cone}) $\Rightarrow$ (\ref{thm:SL-DMM-neg-cone:FSI}):  Let $\alg{B}$ be a totally ordered negatively generated De Morgan monoid that is not a Sugihara monoid.
Then, by Theorem~\ref{thm:FSI-DMMs}, $\alg{B} \cong \alg{S}[\alg{B}']$ for a nontrivial anti-idempotent subalgebra $\alg{B}'$ of $\alg{B}$ and an odd Sugihara monoid $\alg{S}$ (both totally ordered). Suppose, with a view to contradiction, that $\alg{B}'$ is not negatively generated, i.e., $\alg{B}'' \seteq \SgA^{\alg{B}'} B^{-}$ is a proper subalgebra of $\alg{B}'$.
Then, by Theorem~\ref{thm:expandingConstruction}(\ref{thm:expandingConstruction:S}), $\alg{S}[\alg{B}'']$ is a proper subalgebra of $\alg{S}[\alg{B}']$ containing $S[B']^{-}$, contradicting the fact that $\alg{S}[\alg{B}']$ is negatively generated.
So, $\alg{B}'$ is negatively generated, totally ordered, and anti-idempotent, whence $\alg{B}' \cong \refl(\alg{D})$ for some totally ordered $\alg{D} \in \cl{GSM}$, by Theorem~\ref{thm:SLaiDMM-neg-cone}.

(\ref{thm:SL-DMM-neg-cone:FSI}) $\Rightarrow$ (\ref{thm:SL-DMM-neg-cone:equation}): First we show that (\ref{eq:negcone}) holds for every Sugihara monoid, using the fact that $\cl{SM} = \Vop(\alg{Z}^{*})$. For $a \in Z^{*}$, we have $d(a) = a \wedge \neg a = d'(a)$, $\sigma(a) = a = \sigma'(a)$, and
$
(f^2 \vee \neg(f^2)) \to a = e \to a = a.
$
Therefore, the right-hand side of (\ref{eq:negcone}) simplifies to $(a \wedge \neg a) \vee a$, which clearly equals $a$.

Lastly, let $\alg{B} = \alg{S}[\refl(\alg{D})]$ for some totally ordered $\alg{S} \in \cl{OSM}$ and some totally ordered $\alg{D} \in \cl{GSM}$. We have just seen that $\alg{S}$ satisfies (\ref{eq:negcone}). And by Theorem~\ref{thm:SLaiDMM-neg-cone}, the subalgebra $\refl(\alg{D})$ of $\alg{B}$ also satisfies (\ref{eq:negcone}).

Let $a \in B \setminus \refl(D)$, and let $b$ be the right-hand side of (\ref{eq:negcone}) when $x$ is assigned the value of $a$. Recall from Theorems~\ref{thm:Sugihara-factor-odd} and~\ref{thm:FSI-DMMs} that there is a homomorphism from $\alg{B}$ onto $\alg{S}$, whose kernel identifies two distinct elements iff they belong to $\refl(D)$.
So, if $a \neq b$, then since $a \notin \refl(D)$, it follows that $h(a)$ is not $h(b)$, contradicting the fact that $\alg{S}$ satisfies (\ref{eq:negcone}).

The proof of (\ref{thm:SL-DMM-neg-cone:equation})~$\Rightarrow$~(\ref{thm:SL-DMM-neg-cone:neg-cone}) is similiar to its counterpart in
Theorem~\ref{thm:SLaiDMM-neg-cone}.
\end{proof}

\begin{corollary}
\label{cor:SL-DMM-loc-fin}
Let\/ $\cl{K}$ be the class of all negatively generated semilinear De Morgan monoids.
\begin{enumerate}
\item\label{cor:SL-DMM-loc-fin:variety} $\cl{K}$ is a variety and it is axiomatized relative to semilinear De Morgan monoids by \textup{(\ref{eq:negcone}).}
\item \label{cor:SL-DMM-loc-fin:bounds} If $\alg{A} \in \cl{K}$ is totally ordered and $n$-generated, then $|A| \leq 6n + 4$.
\item \label{cor:SL-DMM-loc-fin:loc-fin} $\cl{K}$ is locally finite.
\end{enumerate}
\end{corollary}

\begin{proof}
(\ref{cor:SL-DMM-loc-fin:variety}) follows directly from Theorem~\ref{thm:SL-DMM-neg-cone}.

Let $\alg{A} \in \cl{K}$ be totally ordered an $n$-generated, where $n \in \omega$. If $\alg{A}$ is a Sugihara monoid, then $|A| \leq 2n +2 \leq 6n + 4$ (see Theorem~\ref{thm:FG-SI-SM}). If $\alg{A}$ is not a Sugihara monoid, then $\alg{A} \cong \alg{S}[\alg{A}']$ for an anti-idempotent subalgebra $\alg{A}'$ of $\alg{A}$, and a totally ordered odd Sugihara monoid $\alg{S}$, by Theorem~\ref{thm:FSI-DMMs}. Let us divide the
generators of $\alg{S}[\alg{A}']$ into $X \subseteq A'$ and $Y \subseteq S \setminus \{e^{\alg{S}}\}$, so that when $|X| = p$ and $|Y| = q$, we have $p + q \leq n$. Since $\alg{A}'$ is totally ordered, anti-idempotent and negatively generated, $|A'| \leq 6p +4$, by Corollary~\ref{cor:SlaiDMM-loc-fin}.
Now $\alg{S}$ is generated by $Y$, because if some proper subalgebra $\alg{P}$ of $\alg{S}$ contained $Y$ then, by Theorem~\ref{thm:expandingConstruction}(\ref{thm:expandingConstruction:S}), $\alg{P}[\alg{A}']$ would be a proper subalgebra of $\alg{S}[\alg{A}']$ containing $X \cup Y$, a contradiction.
So, by Theorem~\ref{thm:FG-SI-SM}, $\left|S \setminus \{e^{\alg{S}}\}\right| \leq 2q$. But then $|A| \leq 2q + 6p +4 \leq 6(p+q) +4 \leq 6n +4$, proving (\ref{cor:SL-DMM-loc-fin:bounds}).

Therefore, $\cl{K}$ is locally finite, by Lemma~\ref{fct:FG-function}
(since the SI algebras in $\cl{K}$ are totally ordered).
\end{proof}

Now we can strengthen Theorem~\ref{thm:SLaiDMM-ES} as follows:

\begin{theorem}
\label{thm:SL-DMM-ES}
Let\/ $\cl{K}$ be any variety of negatively generated semilinear De Morgan monoids.
Then\/ $\cl{K}$ has surjective epimorphisms.
\end{theorem}

\begin{proof}
Suppose not. By Theorem~\ref{thm:Campercholi},
there exists $\alg{A} \in \cl{K}_\textup{FSI}$ with a $\cl{K}$-epic proper subalgebra $\alg{B}$.
We proceed to derive a contradiction.

Let $\cl{K}^{SM}$ be the class of all idempotent members of $\cl{K}$. 
Note that $\cl{K}^{SM}$ is a variety of Sugihara monoids, so it has surjective epimorphisms, by Theorem~\ref{thm:SM-Es}. Therefore, $\alg{A}$ is not a Sugihara monoid.
Then, by Theorem~\ref{thm:FSI-DMMs},
$\alg{A} = \alg{S}[\alg{A}']$ for some nontrivial anti-idempotent $\alg{A}' \in \cl{K}$ and some odd Sugihara monoid $\alg{S}$, both totally ordered.

Let $\alg{B}'$ be the subalgebra of $\alg{A}'$ with universe $A' \cap B$.
We show
that $\alg{B}'$ is a $\Vop(\alg{A}')$-epic proper subalgebra of $\alg{A}'$.  This will conclude the proof, as it contradicts the fact that $\Vop(\alg{A}')$ has surjective epimorphisms (by Theorem~\ref{thm:SLaiDMM-ES}).

First, we claim that $\alg{B} = \alg{S}[\alg{B}']$. Evidently $B \subseteq S[B']$. Suppose, with a view to contradiction, that $a \in S \setminus B$. Note that $a \in A$. Let $h \colon \alg{A} \rightarrow \alg{S}$ be the extension, from Theorem~\ref{thm:expandingConstruction}(\ref{thm:expandingConstruction:H}), of the homomorphism which maps $\alg{A}'$ onto the trivial algebra. Then $h(a) \notin h[B]$, by definition of $h$, since $a \in S$. It therefore follows from the surjectivity of $h$ that $h[B]$ is a $\cl{K}^{SM}$-epic proper subalgebra of $\alg{S}$ (since $\alg{B}$ is $\cl{K}$-epic in $\alg{A}$ and compositions of epimorphisms are epimorphisms).
But then $\cl{K}^{SM}$ does not have the ES property, a contradiction. This confirms that $\alg{B} = \alg{S}[\alg{B}']$.

Since $B \subsetneq A = S[A']$, it follows from the claim just proved that $B' \subsetneq A'$, so it remains only to show that $\alg{B}'$ is $\Vop(\alg{A}')$-epic in $\alg{A}'$. Let $h,g \colon \alg{A}' \rightarrow \alg{C}$ be homomorphisms into some $\alg{C} \in \Vop(\alg{A}')_\textup{SI}$ such that $h|_{B'} = g|_{B'}$. By J\'onsson's Theorem,
$\alg{C} \in \Hop \Sop \Puop (\alg{A}')$. By Corollary~\ref{cor:expandingHSPu}, $\alg{S}[\alg{C}] \in \Hop \Sop \Puop (\alg{S}[\alg{A}']) \subseteq \cl{K}$. We extend $h$ and $g$ to homomorphisms $h'$ and $g'$ from $\alg{S}[\alg{A}']$ to $\alg{S}[\alg{C}]$, as in Theorem~\ref{thm:expandingConstruction}(\ref{thm:expandingConstruction:H}). Note that $h'|_{B} = g'|_{B}$, because $B = S[B']$, and $h'|_{S} = g'|_{S}$, by construction. But then $h' = g'$, since $\alg{B}$ is $\cl{K}$-epic in $\alg{A}$. Therefore, $h=g$, so $\alg{B}'$ is $\Vop(\alg{A}')$-epic in $\alg{A}'$, by Lemma~\ref{lem:FSI-image}. 
\end{proof}

\smallskip

\noindent
\textbf{Declarations.}
\begin{itemize}
\item
Funding: See footnote to first page.

\item
Conflicts of interest/Competing interests: The authors declare that there
is no conflict of interest.  The second author is a member of the editorial board of Algebra Universalis.
\end{itemize}

\end{document}